%%%%%%%%%%%%%%%%%%%%%%%%%%%%%%%%%%%%%%%%%%%%%%%%%%%%%%%%%%%%
%%%%  Differential geometry of Cartan connections  %%%%%%%%%
%%%%  Dmitri V\. Alekseevsky and Peter W\. Michor  %%%%%%%%%
%%%%  AMSTeX version 2.1, uses also diag.tex       %%%%%%%%%
%%%%  from Lamstex, 23 pages                       %%%%%%%%%
%%%%%%%%%%%%%%%%%%%%%%%%%%%%%%%%%%%%%%%%%%%%%%%%%%%%%%%%%%%%
% TeX-NMB program applied to the file from 1993.7.16;16:34 on 1993.7.16; 16:34
% TeX-NMB program applied to the file from 1993.7.9;15:51 on 1993.7.9; 15:51
% TeX-NMB program applied to the file from 1993.5.27;15:25 on 1993.5.27; 15:25
% TeX-NMB program applied to the file from 1993.2.18;9:25 on 1993.2.18; 9:26
\input amstex % version 2.1
\input amsppt.sty % version 2.1
%\input diag % from lamstex, for making diagrams: newCD instead of CD 
% macro package for drawing diagrams --
% essentially a part of lamstex.tex
% use exactly following the lamstex manual 
% with '..CD..' replaced by  '..newCD..'
% the original construction $\CD ... \endCD$ remains without any 
% change
%
\catcode`\@=11

\def\East#1#2{\setboxz@h{$\m@th\ssize\;{#1}\;\;$}%
 \setbox@ne\hbox{$\m@th\ssize\;{#2}\;\;$}\setbox\tw@\hbox{$\m@th#2$}%
 \dimen@\minaw@
 \ifdim\wdz@>\dimen@ \dimen@\wdz@ \fi  \ifdim\wd@ne>\dimen@ \dimen@\wd@ne \fi
 \ifdim\wd\tw@>\z@
  \mathrel{\mathop{\hbox to\dimen@{\rightarrowfill}}\limits^{#1}_{#2}}%
 \else
  \mathrel{\mathop{\hbox to\dimen@{\rightarrowfill}}\limits^{#1}}%
 \fi}
\def\West#1#2{\setboxz@h{$\m@th\ssize\;\;{#1}\;$}%
 \setbox@ne\hbox{$\m@th\ssize\;\;{#2}\;$}\setbox\tw@\hbox{$\m@th#2$}%
 \dimen@\minaw@
 \ifdim\wdz@>\dimen@ \dimen@\wdz@ \fi \ifdim\wd@ne>\dimen@ \dimen@\wd@ne \fi
 \ifdim\wd\tw@>\z@
  \mathrel{\mathop{\hbox to\dimen@{\leftarrowfill}}\limits^{#1}_{#2}}%
 \else
  \mathrel{\mathop{\hbox to\dimen@{\leftarrowfill}}\limits^{#1}}%
 \fi}
\font\arrow@i=lams1
\font\arrow@ii=lams2
\font\arrow@iii=lams3
\font\arrow@iv=lams4
\font\arrow@v=lams5
\newbox\zer@
\newdimen\standardcgap
\standardcgap=40\p@
\newdimen\hunit
\hunit=\tw@\p@
\newdimen\standardrgap
\standardrgap=32\p@
\newdimen\vunit
\vunit=1.6\p@
\def\Cgaps#1{\RIfM@
  \standardcgap=#1\standardcgap\relax \hunit=#1\hunit\relax
 \else \nonmatherr@\Cgaps \fi}
\def\Rgaps#1{\RIfM@
  \standardrgap=#1\standardrgap\relax \vunit=#1\vunit\relax
 \else \nonmatherr@\Rgaps \fi}
\newdimen\getdim@
\def\getcgap@#1{\ifcase#1\or\getdim@\z@\else\getdim@\standardcgap\fi}
\def\getrgap@#1{\ifcase#1\getdim@\z@\else\getdim@\standardrgap\fi}
\def\cgaps#1{\RIfM@
 \cgaps@{#1}\edef\getcgap@##1{\i@=##1\relax\the\toks@}\toks@{}\else
 \nonmatherr@\cgaps\fi}
\def\rgaps#1{\RIfM@
 \rgaps@{#1}\edef\getrgap@##1{\i@=##1\relax\the\toks@}\toks@{}\else
 \nonmatherr@\rgaps\fi}
\def\Gaps@@{\gaps@@}
\def\cgaps@#1{\toks@{\ifcase\i@\or\getdim@=\z@}%
 \gaps@@\standardcgap#1;\gaps@@\gaps@@
 \edef\next@{\the\toks@\noexpand\else\noexpand\getdim@\noexpand\standardcgap
  \noexpand\fi}%
 \toks@=\expandafter{\next@}}
\def\rgaps@#1{\toks@{\ifcase\i@\getdim@=\z@}%
 \gaps@@\standardrgap#1;\gaps@@\gaps@@
 \edef\next@{\the\toks@\noexpand\else\noexpand\getdim@\noexpand\standardrgap
  \noexpand\fi}%
 \toks@=\expandafter{\next@}}
\def\gaps@@#1#2;#3{\mgaps@#1#2\mgaps@
 \edef\next@{\the\toks@\noexpand\or\noexpand\getdim@
  \noexpand#1\the\mgapstoks@@}%
 \global\toks@=\expandafter{\next@}%
 \DN@{#3}%
 \ifx\next@\Gaps@@\gdef\next@##1\gaps@@{}\else
  \gdef\next@{\gaps@@#1#3}\fi\next@}
\def\mgaps@#1{\let\mgapsnext@#1\FN@\mgaps@@}
\def\mgaps@@{\ifx\next\space@\DN@. {\FN@\mgaps@@}\else
 \DN@.{\FN@\mgaps@@@}\fi\next@.}
\def\mgaps@@@{\ifx\next\w\let\next@\mgaps@@@@\else
 \let\next@\mgaps@@@@@\fi\next@}
\newtoks\mgapstoks@@
\def\mgaps@@@@@#1\mgaps@{\getdim@\mgapsnext@\getdim@#1\getdim@
 \edef\next@{\noexpand\getdim@\the\getdim@}%
 \mgapstoks@@=\expandafter{\next@}}
\def\mgaps@@@@\w#1#2\mgaps@{\mgaps@@@@@#2\mgaps@
 \setbox\zer@\hbox{$\m@th\hskip15\p@\tsize@#1$}%
 \dimen@\wd\zer@
 \ifdim\dimen@>\getdim@ \getdim@\dimen@ \fi
 \edef\next@{\noexpand\getdim@\the\getdim@}%
 \mgapstoks@@=\expandafter{\next@}}
\def\changewidth#1#2{\setbox\zer@\hbox{$\m@th#2}%
 \hbox to\wd\zer@{\hss$\m@th#1$\hss}}
\atdef@({\FN@\ARROW@}
\def\ARROW@{\ifx\next)\let\next@\OPTIONS@\else
 \DN@{\csname\string @(\endcsname}\fi\next@}
\newif\ifoptions@
\def\OPTIONS@){\ifoptions@\let\next@\relax\else
 \DN@{\options@true\begingroup\optioncodes@}\fi\next@}
\newif\ifN@
\newif\ifE@
\newif\ifNESW@
\newif\ifH@
\newif\ifV@
\newif\ifHshort@
\expandafter\def\csname\string @(\endcsname #1,#2){%
 \ifoptions@\let\next@\endgroup\else\let\next@\relax\fi\next@
 \N@false\E@false\H@false\V@false\Hshort@false
 \ifnum#1>\z@\E@true\fi
 \ifnum#1=\z@\V@true\tX@false\tY@false\a@false\fi
 \ifnum#2>\z@\N@true\fi
 \ifnum#2=\z@\H@true\tX@false\tY@false\a@false\ifshort@\Hshort@true\fi\fi
 \NESW@false
 \ifN@\ifE@\NESW@true\fi\else\ifE@\else\NESW@true\fi\fi
 \arrow@{#1}{#2}%
 \global\options@false
 \global\scount@\z@\global\tcount@\z@\global\arrcount@\z@
 \global\s@false\global\sxdimen@\z@\global\sydimen@\z@
 \global\tX@false\global\tXdimen@i\z@\global\tXdimen@ii\z@
 \global\tY@false\global\tYdimen@i\z@\global\tYdimen@ii\z@
 \global\a@false\global\exacount@\z@
 \global\x@false\global\xdimen@\z@
 \global\X@false\global\Xdimen@\z@
 \global\y@false\global\ydimen@\z@
 \global\Y@false\global\Ydimen@\z@
 \global\p@false\global\pdimen@\z@
 \global\label@ifalse\global\label@iifalse
 \global\dl@ifalse\global\ldimen@i\z@
 \global\dl@iifalse\global\ldimen@ii\z@
 \global\short@false\global\unshort@false}
\newif\iflabel@i
\newif\iflabel@ii
\newcount\scount@
\newcount\tcount@
\newcount\arrcount@
\newif\ifs@
\newdimen\sxdimen@
\newdimen\sydimen@
\newif\iftX@
\newdimen\tXdimen@i
\newdimen\tXdimen@ii
\newif\iftY@
\newdimen\tYdimen@i
\newdimen\tYdimen@ii
\newif\ifa@
\newcount\exacount@
\newif\ifx@
\newdimen\xdimen@
\newif\ifX@
\newdimen\Xdimen@
\newif\ify@
\newdimen\ydimen@
\newif\ifY@
\newdimen\Ydimen@
\newif\ifp@
\newdimen\pdimen@
\newif\ifdl@i
\newif\ifdl@ii
\newdimen\ldimen@i
\newdimen\ldimen@ii
\newif\ifshort@
\newif\ifunshort@
\def\zero@#1{\ifnum\scount@=\z@
 \if#1e\global\scount@\m@ne\else
 \if#1t\global\scount@\tw@\else
 \if#1h\global\scount@\thr@@\else
 \if#1'\global\scount@6 \else
 \if#1`\global\scount@7 \else
 \if#1(\global\scount@8 \else
 \if#1)\global\scount@9 \else
 \if#1s\global\scount@12 \else
 \if#1H\global\scount@13 \else
 \Err@{\Invalid@@ option \string\0}\fi\fi\fi\fi\fi\fi\fi\fi\fi
 \fi}
\def\one@#1{\ifnum\tcount@=\z@
 \if#1e\global\tcount@\m@ne\else
 \if#1h\global\tcount@\tw@\else
 \if#1t\global\tcount@\thr@@\else
 \if#1'\global\tcount@4 \else
 \if#1`\global\tcount@5 \else
 \if#1(\global\tcount@10 \else
 \if#1)\global\tcount@11 \else
 \if#1s\global\tcount@12 \else
 \if#1H\global\tcount@13 \else
 \Err@{\Invalid@@ option \string\1}\fi\fi\fi\fi\fi\fi\fi\fi\fi
 \fi}
\def\a@#1{\ifnum\arrcount@=\z@
 \if#10\global\arrcount@\m@ne\else
 \if#1+\global\arrcount@\@ne\else
 \if#1-\global\arrcount@\tw@\else
 \if#1=\global\arrcount@\thr@@\else
 \Err@{\Invalid@@ option \string\a}\fi\fi\fi\fi
 \fi}
\def\ds@(#1;#2){\ifs@\else
 \global\s@true
 \sxdimen@\hunit \global\sxdimen@#1\sxdimen@\relax
 \sydimen@\vunit \global\sydimen@#2\sydimen@\relax
 \fi}
\def\dtX@(#1;#2){\iftX@\else
 \global\tX@true
 \tXdimen@i\hunit \global\tXdimen@i#1\tXdimen@i\relax
 \tXdimen@ii\vunit \global\tXdimen@ii#2\tXdimen@ii\relax
 \fi}
\def\dtY@(#1;#2){\iftY@\else
 \global\tY@true
 \tYdimen@i\hunit \global\tYdimen@i#1\tYdimen@i\relax
 \tYdimen@ii\vunit \global\tYdimen@ii#2\tYdimen@ii\relax
 \fi}
\def\da@#1{\ifa@\else\global\a@true\global\exacount@#1\relax\fi}
\def\dx@#1{\ifx@\else
 \global\x@true
 \xdimen@\hunit \global\xdimen@#1\xdimen@\relax
 \fi}
\def\dX@#1{\ifX@\else
 \global\X@true
 \Xdimen@\hunit \global\Xdimen@#1\Xdimen@\relax
 \fi}
\def\dy@#1{\ify@\else
 \global\y@true
 \ydimen@\vunit \global\ydimen@#1\ydimen@\relax
 \fi}
\def\dY@#1{\ifY@\else
 \global\Y@true
 \Ydimen@\vunit \global\Ydimen@#1\Ydimen@\relax
 \fi}
\def\p@@#1{\ifp@\else
 \global\p@true
 \pdimen@\hunit \divide\pdimen@\tw@ \global\pdimen@#1\pdimen@\relax
 \fi}
\def\L@#1{\iflabel@i\else
 \global\label@itrue \gdef\label@i{#1}%
 \fi}
\def\l@#1{\iflabel@ii\else
 \global\label@iitrue \gdef\label@ii{#1}%
 \fi}
\def\dL@#1{\ifdl@i\else
 \global\dl@itrue \ldimen@i\hunit \global\ldimen@i#1\ldimen@i\relax
 \fi}
\def\dl@#1{\ifdl@ii\else
 \global\dl@iitrue \ldimen@ii\hunit \global\ldimen@ii#1\ldimen@ii\relax
 \fi}
\def\s@{\ifunshort@\else\global\short@true\fi}
\def\uns@{\ifshort@\else\global\unshort@true\global\short@false\fi}
\def\optioncodes@{\let\0\zero@\let\1\one@\let\a\a@\let\ds\ds@\let\dtX\dtX@
 \let\dtY\dtY@\let\da\da@\let\dx\dx@\let\dX\dX@\let\dY\dY@\let\dy\dy@
 \let\p\p@@\let\L\L@\let\l\l@\let\dL\dL@\let\dl\dl@\let\s\s@\let\uns\uns@}
\def\slopes@{\\161\\152\\143\\134\\255\\126\\357\\238\\349\\45{10}\\56{11}%
 \\11{12}\\65{13}\\54{14}\\43{15}\\32{16}\\53{17}\\21{18}\\52{19}\\31{20}%
 \\41{21}\\51{22}\\61{23}}
\newcount\tan@i
\newcount\tan@ip
\newcount\tan@ii
\newcount\tan@iip
\newdimen\slope@i
\newdimen\slope@ip
\newdimen\slope@ii
\newdimen\slope@iip
\newcount\angcount@
\newcount\extracount@
\def\slope@{{\slope@i=\secondy@ \advance\slope@i-\firsty@
 \ifN@\else\multiply\slope@i\m@ne\fi
 \slope@ii=\secondx@ \advance\slope@ii-\firstx@
 \ifE@\else\multiply\slope@ii\m@ne\fi
 \ifdim\slope@ii<\z@
  \global\tan@i6 \global\tan@ii\@ne \global\angcount@23
 \else
  \dimen@\slope@i \multiply\dimen@6
  \ifdim\dimen@<\slope@ii
   \global\tan@i\@ne \global\tan@ii6 \global\angcount@\@ne
  \else
   \dimen@\slope@ii \multiply\dimen@6
   \ifdim\dimen@<\slope@i
    \global\tan@i6 \global\tan@ii\@ne \global\angcount@23
   \else
    \tan@ip\z@ \tan@iip \@ne
    \def\\##1##2##3{\global\angcount@=##3\relax
     \slope@ip\slope@i \slope@iip\slope@ii
     \multiply\slope@iip##1\relax \multiply\slope@ip##2\relax
     \ifdim\slope@iip<\slope@ip
      \global\tan@ip=##1\relax \global\tan@iip=##2\relax
     \else
      \global\tan@i=##1\relax \global\tan@ii=##2\relax
      \def\\####1####2####3{}%
     \fi}%
    \slopes@
    \slope@i=\secondy@ \advance\slope@i-\firsty@
    \ifN@\else\multiply\slope@i\m@ne\fi
    \multiply\slope@i\tan@ii \multiply\slope@i\tan@iip \multiply\slope@i\tw@
    \count@\tan@i \multiply\count@\tan@iip
    \extracount@\tan@ip \multiply\extracount@\tan@ii
    \advance\count@\extracount@
    \slope@ii=\secondx@ \advance\slope@ii-\firstx@
    \ifE@\else\multiply\slope@ii\m@ne\fi
    \multiply\slope@ii\count@
    \ifdim\slope@i<\slope@ii
     \global\tan@i=\tan@ip \global\tan@ii=\tan@iip
     \global\advance\angcount@\m@ne
    \fi
   \fi
  \fi
 \fi}%
}
\def\slope@a#1{{\def\\##1##2##3{\ifnum##3=#1\global\tan@i=##1\relax
 \global\tan@ii=##2\relax\fi}\slopes@}}
\newcount\i@
\newcount\j@
\newcount\colcount@
\newcount\Colcount@
\newcount\tcolcount@
\newdimen\rowht@
\newdimen\rowdp@
\newcount\rowcount@
\newcount\Rowcount@
\newcount\maxcolrow@
\newtoks\colwidthtoks@
\newtoks\Rowheighttoks@
\newtoks\Rowdepthtoks@
\newtoks\widthtoks@
\newtoks\Widthtoks@
\newtoks\heighttoks@
\newtoks\Heighttoks@
\newtoks\depthtoks@
\newtoks\Depthtoks@
\newif\iffirstnewCDcr@
\def\dotoks@i{%
 \global\widthtoks@=\expandafter{\the\widthtoks@\else\getdim@\z@\fi}%
 \global\heighttoks@=\expandafter{\the\heighttoks@\else\getdim@\z@\fi}%
 \global\depthtoks@=\expandafter{\the\depthtoks@\else\getdim@\z@\fi}}
\def\dotoks@ii{%
 \global\widthtoks@{\ifcase\j@}%
 \global\heighttoks@{\ifcase\j@}%
 \global\depthtoks@{\ifcase\j@}}
\def\prenewCD@#1\endnewCD{\setbox\zer@
 \vbox{%
  \def\arrow@##1##2{{}}%
  \rowcount@\m@ne \colcount@\z@ \Colcount@\z@
  \firstnewCDcr@true \toks@{}%
  \widthtoks@{\ifcase\j@}%
  \Widthtoks@{\ifcase\i@}%
  \heighttoks@{\ifcase\j@}%
  \Heighttoks@{\ifcase\i@}%
  \depthtoks@{\ifcase\j@}%
  \Depthtoks@{\ifcase\i@}%
  \Rowheighttoks@{\ifcase\i@}%
  \Rowdepthtoks@{\ifcase\i@}%
  \Let@
  \everycr{%
   \noalign{%
    \global\advance\rowcount@\@ne
    \ifnum\colcount@<\Colcount@
    \else
     \global\Colcount@=\colcount@ \global\maxcolrow@=\rowcount@
    \fi
    \global\colcount@\z@
    \iffirstnewCDcr@
     \global\firstnewCDcr@false
    \else
     \edef\next@{\the\Rowheighttoks@\noexpand\or\noexpand\getdim@\the\rowht@}%
      \global\Rowheighttoks@=\expandafter{\next@}%
     \edef\next@{\the\Rowdepthtoks@\noexpand\or\noexpand\getdim@\the\rowdp@}%
      \global\Rowdepthtoks@=\expandafter{\next@}%
     \global\rowht@\z@ \global\rowdp@\z@
     \dotoks@i
     \edef\next@{\the\Widthtoks@\noexpand\or\the\widthtoks@}%
      \global\Widthtoks@=\expandafter{\next@}%
     \edef\next@{\the\Heighttoks@\noexpand\or\the\heighttoks@}%
      \global\Heighttoks@=\expandafter{\next@}%
     \edef\next@{\the\Depthtoks@\noexpand\or\the\depthtoks@}%
      \global\Depthtoks@=\expandafter{\next@}%
     \dotoks@ii
    \fi}}%
  \tabskip\z@
  \halign{&\setbox\zer@\hbox{\vrule height10\p@ width\z@ depth\z@
   $\m@th\displaystyle{##}$}\copy\zer@
   \ifdim\ht\zer@>\rowht@ \global\rowht@\ht\zer@ \fi
   \ifdim\dp\zer@>\rowdp@ \global\rowdp@\dp\zer@ \fi
   \global\advance\colcount@\@ne
   \edef\next@{\the\widthtoks@\noexpand\or\noexpand\getdim@\the\wd\zer@}%
    \global\widthtoks@=\expandafter{\next@}%
   \edef\next@{\the\heighttoks@\noexpand\or\noexpand\getdim@\the\ht\zer@}%
    \global\heighttoks@=\expandafter{\next@}%
   \edef\next@{\the\depthtoks@\noexpand\or\noexpand\getdim@\the\dp\zer@}%
    \global\depthtoks@=\expandafter{\next@}%
   \cr#1\crcr}}%
 \Rowcount@=\rowcount@
 \global\Widthtoks@=\expandafter{\the\Widthtoks@\fi\relax}%
 \edef\Width@##1##2{\i@=##1\relax\j@=##2\relax\the\Widthtoks@}%
 \global\Heighttoks@=\expandafter{\the\Heighttoks@\fi\relax}%
 \edef\Height@##1##2{\i@=##1\relax\j@=##2\relax\the\Heighttoks@}%
 \global\Depthtoks@=\expandafter{\the\Depthtoks@\fi\relax}%
 \edef\Depth@##1##2{\i@=##1\relax\j@=##2\relax\the\Depthtoks@}%
 \edef\next@{\the\Rowheighttoks@\noexpand\fi\relax}%
 \global\Rowheighttoks@=\expandafter{\next@}%
 \edef\Rowheight@##1{\i@=##1\relax\the\Rowheighttoks@}%
 \edef\next@{\the\Rowdepthtoks@\noexpand\fi\relax}%
 \global\Rowdepthtoks@=\expandafter{\next@}%
 \edef\Rowdepth@##1{\i@=##1\relax\the\Rowdepthtoks@}%
 \colwidthtoks@{\fi}%
 \setbox\zer@\vbox{%
  \unvbox\zer@
  \count@\rowcount@
  \loop
   \unskip\unpenalty
   \setbox\zer@\lastbox
   \ifnum\count@>\maxcolrow@ \advance\count@\m@ne
   \repeat
  \hbox{%
   \unhbox\zer@
   \count@\z@
   \loop
    \unskip
    \setbox\zer@\lastbox
    \edef\next@{\noexpand\or\noexpand\getdim@\the\wd\zer@\the\colwidthtoks@}%
     \global\colwidthtoks@=\expandafter{\next@}%
    \advance\count@\@ne
    \ifnum\count@<\Colcount@
    \repeat}}%
 \edef\next@{\noexpand\ifcase\noexpand\i@\the\colwidthtoks@}%
  \global\colwidthtoks@=\expandafter{\next@}%
 \edef\Colwidth@##1{\i@=##1\relax\the\colwidthtoks@}%
 \colwidthtoks@{}\Rowheighttoks@{}\Rowdepthtoks@{}\widthtoks@{}%
 \Widthtoks@{}\heighttoks@{}\Heighttoks@{}\depthtoks@{}\Depthtoks@{}%
}
\newcount\xoff@
\newcount\yoff@
\newcount\endcount@
\newcount\rcount@
\newdimen\firstx@
\newdimen\firsty@
\newdimen\secondx@
\newdimen\secondy@
\newdimen\tocenter@
\newdimen\charht@
\newdimen\charwd@
\def\outside@{\Err@{This arrow points outside the \string\newCD}}
\newif\ifsvertex@
\newif\iftvertex@
\def\arrow@#1#2{\xoff@=#1\relax\yoff@=#2\relax
 \count@\rowcount@ \advance\count@-\yoff@
 \ifnum\count@<\@ne \outside@ \else \ifnum\count@>\Rowcount@ \outside@ \fi\fi
 \count@\colcount@ \advance\count@\xoff@
 \ifnum\count@<\@ne \outside@ \else \ifnum\count@>\Colcount@ \outside@\fi\fi
 \tcolcount@\colcount@ \advance\tcolcount@\xoff@
 \Width@\rowcount@\colcount@ \tocenter@=-\getdim@ \divide\tocenter@\tw@
 \ifdim\getdim@=\z@
  \firstx@\z@ \firsty@\mathaxis@ \svertex@true
 \else
  \svertex@false
  \ifHshort@
   \Colwidth@\colcount@
    \ifE@ \firstx@=.5\getdim@ \else \firstx@=-.5\getdim@ \fi
  \else
   \ifE@ \firstx@=\getdim@ \else \firstx@=-\getdim@ \fi
   \divide\firstx@\tw@
  \fi
  \ifE@
   \ifH@ \advance\firstx@\thr@@\p@ \else \advance\firstx@-\thr@@\p@ \fi
  \else
   \ifH@ \advance\firstx@-\thr@@\p@ \else \advance\firstx@\thr@@\p@ \fi
  \fi
  \ifN@
   \Height@\rowcount@\colcount@ \firsty@=\getdim@
   \ifV@ \advance\firsty@\thr@@\p@ \fi
  \else
   \ifV@
    \Depth@\rowcount@\colcount@ \firsty@=-\getdim@
    \advance\firsty@-\thr@@\p@
   \else
    \firsty@\z@
   \fi
  \fi
 \fi
 \ifV@
 \else
  \Colwidth@\colcount@
  \ifE@ \secondx@=\getdim@ \else \secondx@=-\getdim@ \fi
  \divide\secondx@\tw@
  \ifE@ \else \getcgap@\colcount@ \advance\secondx@-\getdim@ \fi
  \endcount@=\colcount@ \advance\endcount@\xoff@
  \count@=\colcount@
  \ifE@
   \advance\count@\@ne
   \loop
    \ifnum\count@<\endcount@
    \Colwidth@\count@ \advance\secondx@\getdim@
    \getcgap@\count@ \advance\secondx@\getdim@
    \advance\count@\@ne
    \repeat
  \else
   \advance\count@\m@ne
   \loop
    \ifnum\count@>\endcount@
    \Colwidth@\count@ \advance\secondx@-\getdim@
    \getcgap@\count@ \advance\secondx@-\getdim@
    \advance\count@\m@ne
    \repeat
  \fi
  \Colwidth@\count@ \divide\getdim@\tw@
  \ifHshort@
  \else
   \ifE@ \advance\secondx@\getdim@ \else \advance\secondx@-\getdim@ \fi
  \fi
  \ifE@ \getcgap@\count@ \advance\secondx@\getdim@ \fi
  \rcount@\rowcount@ \advance\rcount@-\yoff@
  \Width@\rcount@\count@ \divide\getdim@\tw@
  \tvertex@false
  \ifH@\ifdim\getdim@=\z@\tvertex@true\Hshort@false\fi\fi
  \ifHshort@
  \else
   \ifE@ \advance\secondx@-\getdim@ \else \advance\secondx@\getdim@ \fi
  \fi
  \iftvertex@
   \advance\secondx@.4\p@
  \else
   \ifE@ \advance\secondx@-\thr@@\p@ \else \advance\secondx@\thr@@\p@ \fi
  \fi
 \fi
 \ifH@
 \else
  \ifN@
   \Rowheight@\rowcount@ \secondy@\getdim@
  \else
   \Rowdepth@\rowcount@ \secondy@-\getdim@
   \getrgap@\rowcount@ \advance\secondy@-\getdim@
  \fi
  \endcount@=\rowcount@ \advance\endcount@-\yoff@
  \count@=\rowcount@
  \ifN@
   \advance\count@\m@ne
   \loop
    \ifnum\count@>\endcount@
    \Rowheight@\count@ \advance\secondy@\getdim@
    \Rowdepth@\count@ \advance\secondy@\getdim@
    \getrgap@\count@ \advance\secondy@\getdim@
    \advance\count@\m@ne
    \repeat
  \else
   \advance\count@\@ne
   \loop
    \ifnum\count@<\endcount@
    \Rowheight@\count@ \advance\secondy@-\getdim@
    \Rowdepth@\count@ \advance\secondy@-\getdim@
    \getrgap@\count@ \advance\secondy@-\getdim@
    \advance\count@\@ne
    \repeat
  \fi
  \tvertex@false
  \ifV@\Width@\count@\colcount@\ifdim\getdim@=\z@\tvertex@true\fi\fi
  \ifN@
   \getrgap@\count@ \advance\secondy@\getdim@
   \Rowdepth@\count@ \advance\secondy@\getdim@
   \iftvertex@
    \advance\secondy@\mathaxis@
   \else
    \Depth@\count@\tcolcount@ \advance\secondy@-\getdim@
    \advance\secondy@-\thr@@\p@
   \fi
  \else
   \Rowheight@\count@ \advance\secondy@-\getdim@
   \iftvertex@
    \advance\secondy@\mathaxis@
   \else
    \Height@\count@\tcolcount@ \advance\secondy@\getdim@
    \advance\secondy@\thr@@\p@
   \fi
  \fi
 \fi
 \ifV@\else\advance\firstx@\sxdimen@\fi
 \ifH@\else\advance\firsty@\sydimen@\fi
 \iftX@
  \advance\secondy@\tXdimen@ii
  \advance\secondx@\tXdimen@i
  \slope@
 \else
  \iftY@
   \advance\secondy@\tYdimen@ii
   \advance\secondx@\tYdimen@i
   \slope@
   \secondy@=\secondx@ \advance\secondy@-\firstx@
   \ifNESW@ \else \multiply\secondy@\m@ne \fi
   \multiply\secondy@\tan@i \divide\secondy@\tan@ii \advance\secondy@\firsty@
  \else
   \ifa@
    \slope@
    \ifNESW@ \global\advance\angcount@\exacount@ \else
      \global\advance\angcount@-\exacount@ \fi
    \ifnum\angcount@>23 \angcount@23 \fi
    \ifnum\angcount@<\@ne \angcount@\@ne \fi
    \slope@a\angcount@
    \ifY@
     \advance\secondy@\Ydimen@
    \else
     \ifX@
      \advance\secondx@\Xdimen@
      \dimen@\secondx@ \advance\dimen@-\firstx@
      \ifNESW@\else\multiply\dimen@\m@ne\fi
      \multiply\dimen@\tan@i \divide\dimen@\tan@ii
      \advance\dimen@\firsty@ \secondy@=\dimen@
     \fi
    \fi
   \else
    \ifH@\else\ifV@\else\slope@\fi\fi
   \fi
  \fi
 \fi
 \ifH@\else\ifV@\else\ifsvertex@\else
  \dimen@=6\p@ \multiply\dimen@\tan@ii
  \count@=\tan@i \advance\count@\tan@ii \divide\dimen@\count@
  \ifE@ \advance\firstx@\dimen@ \else \advance\firstx@-\dimen@ \fi
  \multiply\dimen@\tan@i \divide\dimen@\tan@ii
  \ifN@ \advance\firsty@\dimen@ \else \advance\firsty@-\dimen@ \fi
 \fi\fi\fi
 \ifp@
  \ifH@\else\ifV@\else
   \getcos@\pdimen@ \advance\firsty@\dimen@ \advance\secondy@\dimen@
   \ifNESW@ \advance\firstx@-\dimen@ii \else \advance\firstx@\dimen@ii \fi
  \fi\fi
 \fi
 \ifH@\else\ifV@\else
  \ifnum\tan@i>\tan@ii
   \charht@=10\p@ \charwd@=10\p@
   \multiply\charwd@\tan@ii \divide\charwd@\tan@i
  \else
   \charwd@=10\p@ \charht@=10\p@
   \divide\charht@\tan@ii \multiply\charht@\tan@i
  \fi
  \ifnum\tcount@=\thr@@
   \ifN@ \advance\secondy@-.3\charht@ \else\advance\secondy@.3\charht@ \fi
  \fi
  \ifnum\scount@=\tw@
   \ifE@ \advance\firstx@.3\charht@ \else \advance\firstx@-.3\charht@ \fi
  \fi
  \ifnum\tcount@=12
   \ifN@ \advance\secondy@-\charht@ \else \advance\secondy@\charht@ \fi
  \fi
  \iftY@
  \else
   \ifa@
    \ifX@
    \else
     \secondx@\secondy@ \advance\secondx@-\firsty@
     \ifNESW@\else\multiply\secondx@\m@ne\fi
     \multiply\secondx@\tan@ii \divide\secondx@\tan@i
     \advance\secondx@\firstx@
    \fi
   \fi
  \fi
 \fi\fi
 \ifH@\harrow@\else\ifV@\varrow@\else\arrow@@\fi\fi}
\newdimen\mathaxis@
\mathaxis@90\p@ \divide\mathaxis@36
\def\harrow@b{\ifE@\hskip\tocenter@\hskip\firstx@\fi}
\def\harrow@bb{\ifE@\hskip\xdimen@\else\hskip\Xdimen@\fi}
\def\harrow@e{\ifE@\else\hskip-\firstx@\hskip-\tocenter@\fi}
\def\harrow@ee{\ifE@\hskip-\Xdimen@\else\hskip-\xdimen@\fi}
\def\harrow@{\dimen@\secondx@\advance\dimen@-\firstx@
 \ifE@ \let\next@\rlap \else  \multiply\dimen@\m@ne \let\next@\llap \fi
 \next@{%
  \harrow@b
  \smash{\raise\pdimen@\hbox to\dimen@
   {\harrow@bb\arrow@ii
    \ifnum\arrcount@=\m@ne \else \ifnum\arrcount@=\thr@@ \else
     \ifE@
      \ifnum\scount@=\m@ne
      \else
       \ifcase\scount@\or\or\char118 \or\char117 \or\or\or\char119 \or
       \char120 \or\char121 \or\char122 \or\or\or\arrow@i\char125 \or
       \char117 \hskip\thr@@\p@\char117 \hskip-\thr@@\p@\fi
      \fi
     \else
      \ifnum\tcount@=\m@ne
      \else
       \ifcase\tcount@\char117 \or\or\char117 \or\char118 \or\char119 \or
       \char120\or\or\or\or\or\char121 \or\char122 \or\arrow@i\char125
       \or\char117 \hskip\thr@@\p@\char117 \hskip-\thr@@\p@\fi
      \fi
     \fi
    \fi\fi
    \dimen@\mathaxis@ \advance\dimen@.2\p@
    \dimen@ii\mathaxis@ \advance\dimen@ii-.2\p@
    \ifnum\arrcount@=\m@ne
     \let\leads@\null
    \else
     \ifcase\arrcount@
      \def\leads@{\hrule height\dimen@ depth-\dimen@ii}\or
      \def\leads@{\hrule height\dimen@ depth-\dimen@ii}\or
      \def\leads@{\hbox to10\p@{%
       \leaders\hrule height\dimen@ depth-\dimen@ii\hfil
       \hfil
      \leaders\hrule height\dimen@ depth-\dimen@ii\hskip\z@ plus2fil\relax
       \hfil
       \leaders\hrule height\dimen@ depth-\dimen@ii\hfil}}\or
     \def\leads@{\hbox{\hbox to10\p@{\dimen@\mathaxis@ \advance\dimen@1.2\p@
       \dimen@ii\dimen@ \advance\dimen@ii-.4\p@
       \leaders\hrule height\dimen@ depth-\dimen@ii\hfil}%
       \kern-10\p@
       \hbox to10\p@{\dimen@\mathaxis@ \advance\dimen@-1.2\p@
       \dimen@ii\dimen@ \advance\dimen@ii-.4\p@
       \leaders\hrule height\dimen@ depth-\dimen@ii\hfil}}}\fi
    \fi
    \cleaders\leads@\hfil
    \ifnum\arrcount@=\m@ne\else\ifnum\arrcount@=\thr@@\else
     \arrow@i
     \ifE@
      \ifnum\tcount@=\m@ne
      \else
       \ifcase\tcount@\char119 \or\or\char119 \or\char120 \or\char121 \or
       \char122 \or \or\or\or\or\char123\or\char124 \or
       \char125 \or\char119 \hskip-\thr@@\p@\char119 \hskip\thr@@\p@\fi
      \fi
     \else
      \ifcase\scount@\or\or\char120 \or\char119 \or\or\or\char121 \or\char122
      \or\char123 \or\char124 \or\or\or\char125 \or
      \char119 \hskip-\thr@@\p@\char119 \hskip\thr@@\p@\fi
     \fi
    \fi\fi
    \harrow@ee}}%
  \harrow@e}%
 \iflabel@i
  \dimen@ii\z@ \setbox\zer@\hbox{$\m@th\tsize@@\label@i$}%
  \ifnum\arrcount@=\m@ne
  \else
   \advance\dimen@ii\mathaxis@
   \advance\dimen@ii\dp\zer@ \advance\dimen@ii\tw@\p@
   \ifnum\arrcount@=\thr@@ \advance\dimen@ii\tw@\p@ \fi
  \fi
  \advance\dimen@ii\pdimen@
  \next@{\harrow@b\smash{\raise\dimen@ii\hbox to\dimen@
   {\harrow@bb\hskip\tw@\ldimen@i\hfil\box\zer@\hfil\harrow@ee}}\harrow@e}%
 \fi
 \iflabel@ii
  \ifnum\arrcount@=\m@ne
  \else
   \setbox\zer@\hbox{$\m@th\tsize@\label@ii$}%
   \dimen@ii-\ht\zer@ \advance\dimen@ii-\tw@\p@
   \ifnum\arrcount@=\thr@@ \advance\dimen@ii-\tw@\p@ \fi
   \advance\dimen@ii\mathaxis@ \advance\dimen@ii\pdimen@
   \next@{\harrow@b\smash{\raise\dimen@ii\hbox to\dimen@
    {\harrow@bb\hskip\tw@\ldimen@ii\hfil\box\zer@\hfil\harrow@ee}}\harrow@e}%
  \fi
 \fi}
\let\tsize@\tsize
\def\tsizenewCDlabels{\let\tsize@\tsize}
\def\ssizenewCDlabels{\let\tsize@\ssize}
\def\tsize@@{\ifnum\arrcount@=\m@ne\else\tsize@\fi}
\def\varrow@{\dimen@\secondy@ \advance\dimen@-\firsty@
 \ifN@ \else \multiply\dimen@\m@ne \fi
 \setbox\zer@\vbox to\dimen@
  {\ifN@ \vskip-\Ydimen@ \else \vskip\ydimen@ \fi
   \ifnum\arrcount@=\m@ne\else\ifnum\arrcount@=\thr@@\else
    \hbox{\arrow@iii
     \ifN@
      \ifnum\tcount@=\m@ne
      \else
       \ifcase\tcount@\char117 \or\or\char117 \or\char118 \or\char119 \or
       \char120 \or\or\or\or\or\char121 \or\char122 \or\char123 \or
       \vbox{\hbox{\char117 }\nointerlineskip\vskip\thr@@\p@
       \hbox{\char117 }\vskip-\thr@@\p@}\fi
      \fi
     \else
      \ifcase\scount@\or\or\char118 \or\char117 \or\or\or\char119 \or
      \char120 \or\char121 \or\char122 \or\or\or\char123 \or
      \vbox{\hbox{\char117 }\nointerlineskip\vskip\thr@@\p@
      \hbox{\char117 }\vskip-\thr@@\p@}\fi
     \fi}%
    \nointerlineskip
   \fi\fi
   \ifnum\arrcount@=\m@ne
    \let\leads@\null
   \else
    \ifcase\arrcount@\let\leads@\vrule\or\let\leads@\vrule\or
    \def\leads@{\vbox to10\p@{%
     \hrule height 1.67\p@ depth\z@ width.4\p@
     \vfil
     \hrule height 3.33\p@ depth\z@ width.4\p@
     \vfil
     \hrule height 1.67\p@ depth\z@ width.4\p@}}\or
    \def\leads@{\hbox{\vrule height\p@\hskip\tw@\p@\vrule}}\fi
   \fi
  \cleaders\leads@\vfill\nointerlineskip
   \ifnum\arrcount@=\m@ne\else\ifnum\arrcount@=\thr@@\else
    \hbox{\arrow@iv
     \ifN@
      \ifcase\scount@\or\or\char118 \or\char117 \or\or\or\char119 \or
      \char120 \or\char121 \or\char122 \or\or\or\arrow@iii\char123 \or
      \vbox{\hbox{\char117 }\nointerlineskip\vskip-\thr@@\p@
      \hbox{\char117 }\vskip\thr@@\p@}\fi
     \else
      \ifnum\tcount@=\m@ne
      \else
       \ifcase\tcount@\char117 \or\or\char117 \or\char118 \or\char119 \or
       \char120 \or\or\or\or\or\char121 \or\char122 \or\arrow@iii\char123 \or
       \vbox{\hbox{\char117 }\nointerlineskip\vskip-\thr@@\p@
       \hbox{\char117 }\vskip\thr@@\p@}\fi
      \fi
     \fi}%
   \fi\fi
   \ifN@\vskip\ydimen@\else\vskip-\Ydimen@\fi}%
 \ifN@
  \dimen@ii\firsty@
 \else
  \dimen@ii-\firsty@ \advance\dimen@ii\ht\zer@ \multiply\dimen@ii\m@ne
 \fi
 \rlap{\smash{\hskip\tocenter@ \hskip\pdimen@ \raise\dimen@ii \box\zer@}}%
 \iflabel@i
  \setbox\zer@\vbox to\dimen@{\vfil
   \hbox{$\m@th\tsize@@\label@i$}\vskip\tw@\ldimen@i\vfil}%
  \rlap{\smash{\hskip\tocenter@ \hskip\pdimen@
  \ifnum\arrcount@=\m@ne \let\next@\relax \else \let\next@\llap \fi
  \next@{\raise\dimen@ii\hbox{\ifnum\arrcount@=\m@ne \hskip-.5\wd\zer@ \fi
   \box\zer@ \ifnum\arrcount@=\m@ne \else \hskip\tw@\p@ \fi}}}}%
 \fi
 \iflabel@ii
  \ifnum\arrcount@=\m@ne
  \else
   \setbox\zer@\vbox to\dimen@{\vfil
    \hbox{$\m@th\tsize@\label@ii$}\vskip\tw@\ldimen@ii\vfil}%
   \rlap{\smash{\hskip\tocenter@ \hskip\pdimen@
   \rlap{\raise\dimen@ii\hbox{\ifnum\arrcount@=\thr@@ \hskip4.5\p@ \else
    \hskip2.5\p@ \fi\box\zer@}}}}%
  \fi
 \fi
}
\newdimen\goal@
\newdimen\shifted@
\newcount\Tcount@
\newcount\Scount@
\newbox\shaft@
\newcount\slcount@
\def\getcos@#1{%
 \ifnum\tan@i<\tan@ii
  \dimen@#1%
  \ifnum\slcount@<8 \count@9 \else \ifnum\slcount@<12 \count@8 \else
   \count@7 \fi\fi
  \multiply\dimen@\count@ \divide\dimen@10
  \dimen@ii\dimen@ \multiply\dimen@ii\tan@i \divide\dimen@ii\tan@ii
 \else
  \dimen@ii#1%
  \count@-\slcount@ \advance\count@24
  \ifnum\count@<8 \count@9 \else \ifnum\count@<12 \count@8
   \else\count@7 \fi\fi
  \multiply\dimen@ii\count@ \divide\dimen@ii10
  \dimen@\dimen@ii \multiply\dimen@\tan@ii \divide\dimen@\tan@i
 \fi}
\newdimen\adjust@
\def\Nnext@{\ifN@\let\next@\raise\else\let\next@\lower\fi}
\def\arrow@@{\slcount@\angcount@
 \ifNESW@
  \ifnum\angcount@<10
   \let\arrowfont@=\arrow@i \advance\angcount@\m@ne \multiply\angcount@13
  \else
   \ifnum\angcount@<19
    \let\arrowfont@=\arrow@ii \advance\angcount@-10 \multiply\angcount@13
   \else
    \let\arrowfont@=\arrow@iii \advance\angcount@-19 \multiply\angcount@13
  \fi\fi
  \Tcount@\angcount@
 \else
  \ifnum\angcount@<5
   \let\arrowfont@=\arrow@iii \advance\angcount@\m@ne \multiply\angcount@13
   \advance\angcount@65
  \else
   \ifnum\angcount@<14
    \let\arrowfont@=\arrow@iv \advance\angcount@-5 \multiply\angcount@13
   \else
    \ifnum\angcount@<23
     \let\arrowfont@=\arrow@v \advance\angcount@-14 \multiply\angcount@13
    \else
     \let\arrowfont@=\arrow@i \angcount@=117
  \fi\fi\fi
  \ifnum\angcount@=117 \Tcount@=115 \else\Tcount@\angcount@ \fi
 \fi
 \Scount@\Tcount@
 \ifE@
  \ifnum\tcount@=\z@ \advance\Tcount@\tw@ \else\ifnum\tcount@=13
   \advance\Tcount@\tw@ \else \advance\Tcount@\tcount@ \fi\fi
  \ifnum\scount@=\z@ \else \ifnum\scount@=13 \advance\Scount@\thr@@ \else
   \advance\Scount@\scount@ \fi\fi
 \else
  \ifcase\tcount@\advance\Tcount@\thr@@\or\or\advance\Tcount@\thr@@\or
  \advance\Tcount@\tw@\or\advance\Tcount@6 \or\advance\Tcount@7
  \or\or\or\or\or \advance\Tcount@8 \or\advance\Tcount@9 \or
  \advance\Tcount@12 \or\advance\Tcount@\thr@@\fi
  \ifcase\scount@\or\or\advance\Scount@\thr@@\or\advance\Scount@\tw@\or
  \or\or\advance\Scount@4 \or\advance\Scount@5 \or\advance\Scount@10
  \or\advance\Scount@11 \or\or\or\advance\Scount@12 \or\advance
  \Scount@\tw@\fi
 \fi
 \ifcase\arrcount@\or\or\advance\angcount@\@ne\else\fi
 \ifN@ \shifted@=\firsty@ \else\shifted@=-\firsty@ \fi
 \ifE@ \else\advance\shifted@\charht@ \fi
 \goal@=\secondy@ \advance\goal@-\firsty@
 \ifN@\else\multiply\goal@\m@ne\fi
 \setbox\shaft@\hbox{\arrowfont@\char\angcount@}%
 \ifnum\arrcount@=\thr@@
  \getcos@{1.5\p@}%
  \setbox\shaft@\hbox to\wd\shaft@{\arrowfont@
   \rlap{\hskip\dimen@ii
    \smash{\ifNESW@\let\next@\lower\else\let\next@\raise\fi
     \next@\dimen@\hbox{\arrowfont@\char\angcount@}}}%
   \rlap{\hskip-\dimen@ii
    \smash{\ifNESW@\let\next@\raise\else\let\next@\lower\fi
      \next@\dimen@\hbox{\arrowfont@\char\angcount@}}}\hfil}%
 \fi
 \rlap{\smash{\hskip\tocenter@\hskip\firstx@
  \ifnum\arrcount@=\m@ne
  \else
   \ifnum\arrcount@=\thr@@
   \else
    \ifnum\scount@=\m@ne
    \else
     \ifnum\scount@=\z@
     \else
      \setbox\zer@\hbox{\ifnum\angcount@=117 \arrow@v\else\arrowfont@\fi
       \char\Scount@}%
      \ifNESW@
       \ifnum\scount@=\tw@
        \dimen@=\shifted@ \advance\dimen@-\charht@
        \ifN@\hskip-\wd\zer@\fi
        \Nnext@
        \next@\dimen@\copy\zer@
        \ifN@\else\hskip-\wd\zer@\fi
       \else
        \Nnext@
        \ifN@\else\hskip-\wd\zer@\fi
        \next@\shifted@\copy\zer@
        \ifN@\hskip-\wd\zer@\fi
       \fi
       \ifnum\scount@=12
        \advance\shifted@\charht@ \advance\goal@-\charht@
        \ifN@ \hskip\wd\zer@ \else \hskip-\wd\zer@ \fi
       \fi
       \ifnum\scount@=13
        \getcos@{\thr@@\p@}%
        \ifN@ \hskip\dimen@ \else \hskip-\wd\zer@ \hskip-\dimen@ \fi
        \adjust@\shifted@ \advance\adjust@\dimen@ii
        \Nnext@
        \next@\adjust@\copy\zer@
        \ifN@ \hskip-\dimen@ \hskip-\wd\zer@ \else \hskip\dimen@ \fi
       \fi
      \else
       \ifN@\hskip-\wd\zer@\fi
       \ifnum\scount@=\tw@
        \ifN@ \hskip\wd\zer@ \else \hskip-\wd\zer@ \fi
        \dimen@=\shifted@ \advance\dimen@-\charht@
        \Nnext@
        \next@\dimen@\copy\zer@
        \ifN@\hskip-\wd\zer@\fi
       \else
        \Nnext@
        \next@\shifted@\copy\zer@
        \ifN@\else\hskip-\wd\zer@\fi
       \fi
       \ifnum\scount@=12
        \advance\shifted@\charht@ \advance\goal@-\charht@
        \ifN@ \hskip-\wd\zer@ \else \hskip\wd\zer@ \fi
       \fi
       \ifnum\scount@=13
        \getcos@{\thr@@\p@}%
        \ifN@ \hskip-\wd\zer@ \hskip-\dimen@ \else \hskip\dimen@ \fi
        \adjust@\shifted@ \advance\adjust@\dimen@ii
        \Nnext@
        \next@\adjust@\copy\zer@
        \ifN@ \hskip\dimen@ \else \hskip-\dimen@ \hskip-\wd\zer@ \fi
       \fi	
      \fi
  \fi\fi\fi\fi
  \ifnum\arrcount@=\m@ne
  \else
   \loop
    \ifdim\goal@>\charht@
    \ifE@\else\hskip-\charwd@\fi
    \Nnext@
    \next@\shifted@\copy\shaft@
    \ifE@\else\hskip-\charwd@\fi
    \advance\shifted@\charht@ \advance\goal@ -\charht@
    \repeat
   \ifdim\goal@>\z@
    \dimen@=\charht@ \advance\dimen@-\goal@
    \divide\dimen@\tan@i \multiply\dimen@\tan@ii
    \ifE@ \hskip-\dimen@ \else \hskip-\charwd@ \hskip\dimen@ \fi
    \adjust@=\shifted@ \advance\adjust@-\charht@ \advance\adjust@\goal@
    \Nnext@
    \next@\adjust@\copy\shaft@
    \ifE@ \else \hskip-\charwd@ \fi
   \else
    \adjust@=\shifted@ \advance\adjust@-\charht@
   \fi
  \fi
  \ifnum\arrcount@=\m@ne
  \else
   \ifnum\arrcount@=\thr@@
   \else
    \ifnum\tcount@=\m@ne
    \else
     \setbox\zer@
      \hbox{\ifnum\angcount@=117 \arrow@v\else\arrowfont@\fi\char\Tcount@}%
     \ifnum\tcount@=\thr@@
      \advance\adjust@\charht@
      \ifE@\else\ifN@\hskip-\charwd@\else\hskip-\wd\zer@\fi\fi
     \else
      \ifnum\tcount@=12
       \advance\adjust@\charht@
       \ifE@\else\ifN@\hskip-\charwd@\else\hskip-\wd\zer@\fi\fi
      \else
       \ifE@\hskip-\wd\zer@\fi
     \fi\fi
     \Nnext@
     \next@\adjust@\copy\zer@
     \ifnum\tcount@=13
      \hskip-\wd\zer@
      \getcos@{\thr@@\p@}%
      \ifE@\hskip-\dimen@ \else\hskip\dimen@ \fi
      \advance\adjust@-\dimen@ii
      \Nnext@
      \next@\adjust@\box\zer@
     \fi
  \fi\fi\fi}}%
 \iflabel@i
  \rlap{\hskip\tocenter@
  \dimen@\firstx@ \advance\dimen@\secondx@ \divide\dimen@\tw@
  \advance\dimen@\ldimen@i
  \dimen@ii\firsty@ \advance\dimen@ii\secondy@ \divide\dimen@ii\tw@
  \multiply\ldimen@i\tan@i \divide\ldimen@i\tan@ii
  \ifNESW@ \advance\dimen@ii\ldimen@i \else \advance\dimen@ii-\ldimen@i \fi
  \setbox\zer@\hbox{\ifNESW@\else\ifnum\arrcount@=\thr@@\hskip4\p@\else
   \hskip\tw@\p@\fi\fi
   $\m@th\tsize@@\label@i$\ifNESW@\ifnum\arrcount@=\thr@@\hskip4\p@\else
   \hskip\tw@\p@\fi\fi}%
  \ifnum\arrcount@=\m@ne
   \ifNESW@ \advance\dimen@.5\wd\zer@ \advance\dimen@\p@ \else
    \advance\dimen@-.5\wd\zer@ \advance\dimen@-\p@ \fi
   \advance\dimen@ii-.5\ht\zer@
  \else
   \advance\dimen@ii\dp\zer@
   \ifnum\slcount@<6 \advance\dimen@ii\tw@\p@ \fi
  \fi
  \hskip\dimen@
  \ifNESW@ \let\next@\llap \else\let\next@\rlap \fi
  \next@{\smash{\raise\dimen@ii\box\zer@}}}%
 \fi
 \iflabel@ii
  \ifnum\arrcount@=\m@ne
  \else
   \rlap{\hskip\tocenter@
   \dimen@\firstx@ \advance\dimen@\secondx@ \divide\dimen@\tw@
   \ifNESW@ \advance\dimen@\ldimen@ii \else \advance\dimen@-\ldimen@ii \fi
   \dimen@ii\firsty@ \advance\dimen@ii\secondy@ \divide\dimen@ii\tw@
   \multiply\ldimen@ii\tan@i \divide\ldimen@ii\tan@ii
   \advance\dimen@ii\ldimen@ii
   \setbox\zer@\hbox{\ifNESW@\ifnum\arrcount@=\thr@@\hskip4\p@\else
    \hskip\tw@\p@\fi\fi
    $\m@th\tsize@\label@ii$\ifNESW@\else\ifnum\arrcount@=\thr@@\hskip4\p@
    \else\hskip\tw@\p@\fi\fi}%
   \advance\dimen@ii-\ht\zer@
   \ifnum\slcount@<9 \advance\dimen@ii-\thr@@\p@ \fi
   \ifNESW@ \let\next@\rlap \else \let\next@\llap \fi
   \hskip\dimen@\next@{\smash{\raise\dimen@ii\box\zer@}}}%
  \fi
 \fi
}
\def\outnewCD@#1{\def#1{\Err@{\string#1 must not be used within \string\newCD}}}
\newskip\prenewCDskip@
\newskip\postnewCDskip@
\prenewCDskip@\z@
\postnewCDskip@\z@
\def\prenewCDspace#1{\RIfMIfI@
 \onlydmatherr@\prenewCDspace\else\advance\prenewCDskip@#1\relax\fi\else
 \onlydmatherr@\prenewCDspace\fi}
\def\postnewCDspace#1{\RIfMIfI@
 \onlydmatherr@\postnewCDspace\else\advance\postnewCDskip@#1\relax\fi\else
 \onlydmatherr@\postnewCDspace\fi}
\def\predisplayspace#1{\RIfMIfI@
 \onlydmatherr@\predisplayspace\else
 \advance\abovedisplayskip#1\relax
 \advance\abovedisplayshortskip#1\relax\fi
 \else\onlydmatherr@\prenewCDspace\fi}
\def\postdisplayspace#1{\RIfMIfI@
 \onlydmatherr@\postdisplayspace\else
 \advance\belowdisplayskip#1\relax
 \advance\belowdisplayshortskip#1\relax\fi
 \else\onlydmatherr@\postdisplayspace\fi}
\def\PrenewCDSpace#1{\global\prenewCDskip@#1\relax}
\def\PostnewCDSpace#1{\global\postnewCDskip@#1\relax}
\def\newCD#1\endnewCD{%
 \outnewCD@\cgaps\outnewCD@\rgaps\outnewCD@\Cgaps\outnewCD@\Rgaps
 \prenewCD@#1\endnewCD
 \advance\abovedisplayskip\prenewCDskip@
 \advance\abovedisplayshortskip\prenewCDskip@
 \advance\belowdisplayskip\postnewCDskip@
 \advance\belowdisplayshortskip\postnewCDskip@
 \vcenter{\vskip\prenewCDskip@ \Let@ \colcount@\@ne \rowcount@\z@
  \everycr{%
   \noalign{%
    \ifnum\rowcount@=\Rowcount@
    \else
     \global\nointerlineskip
     \getrgap@\rowcount@ \vskip\getdim@
     \global\advance\rowcount@\@ne \global\colcount@\@ne
    \fi}}%
  \tabskip\z@
  \halign{&\global\xoff@\z@ \global\yoff@\z@
   \getcgap@\colcount@ \hskip\getdim@
   \hfil\vrule height10\p@ width\z@ depth\z@
   $\m@th\displaystyle{##}$\hfil
   \global\advance\colcount@\@ne\cr
   #1\crcr}\vskip\postnewCDskip@}%
 \prenewCDskip@\z@\postnewCDskip@\z@
 \def\getcgap@##1{\ifcase##1\or\getdim@\z@\else\getdim@\standardcgap\fi}%
 \def\getrgap@##1{\ifcase##1\getdim@\z@\else\getdim@\standardrgap\fi}%
 \let\Width@\relax\let\Height@\relax\let\Depth@\relax\let\Rowheight@\relax
 \let\Rowdepth@\relax\let\Colwdith@\relax
}
\catcode`\@=\active
%\endinput
%end of diag.tex
\hsize 30pc
\vsize 47pc
\def\nmb#1#2{#2}         % used for renumbering, TeX should ignore.
\def\cit#1#2{\ifx#1!\cite{#2}\else#2\fi}
\def\totoc{}             %= to table of content, invoked by kms-book.sty
\def\idx{}               % for producing index, invoked by kms-book.sty
\def\ign#1{}             %=ignore, invisible entry for the index only

\redefine\o{\circ}
\define\X{\frak X}

\define\ga{\gamma}
\define\de{\delta}

\define\ze{\zeta}

\define\th{\theta}

\define\ka{\kappa}
\define\la{\lambda}
\define\rh{\rho}
\define\si{\sigma}
\define\ta{\tau}
\define\ph{\varphi}
\define\ch{\chi}
\define\ps{\psi}
\define\om{\omega}
\define\Ga{\Gamma}

\define\La{\Lambda}

\define\Ph{\Phi}
\define\Ps{\Psi}
\define\Om{\Omega}
\redefine\i{^{-1}}
\define\row#1#2#3{#1_{#2},\ldots,#1_{#3}}
\define\x{\times}

\define\Fl{\operatorname{Fl}}

\define\Ad{\operatorname{Ad}}
\define\ad{\operatorname{ad}}
\redefine\L{{\Cal L}}
\define\ddt{\left.\tfrac \partial{\partial t}\right\vert_0}
\define\g{{\frak g}}
\define\h{{\frak h}}
\define\pr{{\operatorname{pr}}}
\def\today{\ifcase\month\or
 January\or February\or March\or April\or May\or June\or
 July\or August\or September\or October\or November\or December\fi
 \space\number\day, \number\year}
\topmatter
\title  Differential geometry of Cartan connections
\endtitle
\author Dmitri V\. Alekseevsky \\
Peter W\. Michor  \endauthor
\affil
Erwin Schr\"odinger International Institute of Mathematical Physics, 
Wien, Austria
\endaffil
\address 	 D\. V\. Alekseevsky: 
gen. Antonova 2 kv 99, 117279 Moscow B-279, Russia
\endaddress
\address
P\. W\. Michor: Institut f\"ur Mathematik, Universit\"at Wien,
Strudlhofgasse 4, A-1090 Wien, Austria; and 
Erwin Schr\"odinger International Institute of Mathematical Physics, 
Pasteurgasse 6/7, A-1090 Wien, Austria
\endaddress
\email michor\@esi.ac.at \endemail
%\date {\today} \enddate
\thanks Supported by Project P 7724 PHY
of `Fonds zur F\"orderung der wissenschaftlichen Forschung'.
\endthanks
%\abstract For a more general notion of Cartan connection we define 
%characteristic classes, we investigate their relation to usual 
%characteristic classes, and we $\ldots$.
%\endabstract
\endtopmatter

%\input amspptb.sty
%\userunningheads
%\def\leftheadtext{\smc D\. V\. Alekseevski, P\. W\. Michor}
%\def\rightheadtext{\smc Differential Geometry of Cartan connections}
%\def\bottremark{\today\hfill}

\document

\heading Table of contents \endheading
%\input \jobname.toc
%\loadtoc
%\loadindex
\noindent 1. Introduction \leaders \hbox to 1em{\hss .\hss }\hfill {\eightrm 1}\par 
\noindent 2. Cartan connections and generalized Cartan connections \leaders \hbox to 1em{\hss .\hss }\hfill {\eightrm 3}\par 
\noindent 3. The relation between principal Cartan connections and principal connections \leaders \hbox to 1em{\hss .\hss }\hfill {\eightrm 7}\par 
\noindent 4. Flat Cartan connections \leaders \hbox to 1em{\hss .\hss }\hfill {\eightrm 12}\par 
\noindent 5. Flat Cartan connections associated with a flat $G$-structure \leaders \hbox to 1em{\hss .\hss }\hfill {\eightrm 16}\par 
\noindent 6. The canonical Cartan connection for a $G$-structure of first or second order \leaders \hbox to 1em{\hss .\hss }\hfill {\eightrm 20}\par 

\head\totoc\nmb0{1}. Introduction \endhead

In this article a general theory of Cartan connections is developed 
and some applications are indicated. The starting idea is to consider 
a Cartan connection as a deformation of a local Lie group structure 
on the manifold, i\.e\. a 1-form $\la$ with values in a Lie algebra $\h$ 
which is non degenerate and satisfies the Maurer-Cartan equation. 
Such a Maurer-Cartan form $\la$ may be considered as a flat Cartan 
connection. Many notions and results of the geometry of group 
manifolds are still valid in this more general setting. 

More precisely, for a Lie subalgebra $\g$ of $\h$ we define a Cartan 
connection of type $\h/\g$ on a manifold $P$ of dimension $n=\dim\h$ 
as a $\h$-valued 1-form $\ka:TP\to \h$ which defines an isomorphism 
$\ka_x:T_xP\to \h$ for any $x\in P$ and such that 
$$[\ze_X,\ze_Y]=\ze_{[X,Y]}$$
holds for $X\in \h$ and $Y\in\g$, where the linear mapping 
$\ze:\h\to\X(P)$ from $\h$ into the Lie algebra $\X(P)$ of vector 
fields on $P$ is given by $\ze_X(x)=\ka_x\i(X)$.
If $\g=\h$ then $\ze$ defines a free transitive action of the Lie 
algebra $\h$ on the manifold $P$ in the sense of \cit!{5} and $\ka$ 
is the Maurer-Cartan form of the associated structure of the local 
Lie group structure on $P$. In the general case, when $\g\ne\h$, we 
only have a free action $\ze|\g$ of the Lie algebra $\g$ on $P$. 
So we may think of the Cartan connection $\ka$ as a deformed 
Maurer-Cartan form, where the deformation is breaking the 
symmetry from $\h$ to $\g$. If the action of $\g$ on $P$ can be 
integrated to a free action of a corresponding Lie group $G$ on $P$ with 
smooth orbit space $M=P/G$, the notion of Cartan connection reduces 
to the well known notion of a Cartan connection on the principal 
bundle $p:P\to M$.

In \nmb!{2.3} and \nmb!{2.4} we describe two situations when a Cartan 
connection arises naturally. First under a reduction of a principal 
bundle $p:Q\to M$ with a principal connection to a principal 
subbundle $p:P\to M$. Second when a $G$-structure with a connection 
is given: more precisely, it the Lie algebra admits a reductive 
decomposition $\h=\g\oplus\frak m$ we may identify a Cartan 
connection of type $\h/\g$ on a principal $G$-bundle $p:P\to M$ with 
a $G$-structure on $M$ together with a principal connection in 
$p:P\to M$.

Dropping the condition that the 1-form $\ka$ is non-degenerate we 
come to the notion of \idx{\it generalized Cartan connection}. It is 
closely related with the the notion of a principal connection form on 
a $\g$-manifold, defined in \cit!{5}, see \nmb!{2.6}. In the end of 
section \nmb!{2} we define for an arbitrary generalized Cartan 
connection $\ka$ such notions as the curvature 2-form
$$K=d\ka+\tfrac12[\ka,\ka]^\wedge ,$$
the Bianchi identity
$$dK + [\ka,K]^\wedge =0,$$
the covariant exterior derivative 
$$d_\ka:\Om^p_{\text{hor}}(P;W)^\g\to\Om_{\text{hor}}^{p+1}(M;V)^\g, 
\quad d_\ka(\Ps)= d\Ps + \rh^\wedge(\ka)\Ps,$$
where $\Om_{\text{hor}}^p(M;W)^\g$ is the space of horizontal 
$\g$-equivariant $p$-forms with values in the $\g$-module defined by 
a representation $\rh:\h\to\g\frak l(W)$.

In \nmb!{2.9} we associate with a generalized Cartan connection $\ka$ 
of type $\h/\g$ the Chern-Weil homomorphism
$$\ga:S(\h^*)^\h\to \Om_{\text{hor}}(P)^\g$$
of the algebra of $\h$-invariant polynomials on $\h$ into the algebra 
of $\g$-invariant closed horizontal differential forms on $P$ and 
prove that the characteristic cohomology class $[\ga(f)]$ does not 
depend on the particular choice of the generalized Cartan connection.

In section \nmb!{3} we study relations between principal Cartan 
connections of a principal $G$-bundle $p:P\to M$ and principal 
connections on the $H$-bundle $p:P[H]=P\x_GH\to M$, where $H\supset G$ is a 
Lie group associated to $\h$. We also establish a 
canonical linear isomorphism 
$$\Om_{\text{hor}}(P;W)^G \to \Om_{\text{hor}}(P[H];W)^H$$
between the respective spaces of equivariant horizontal forms with 
values in a representation space $W$ of $H$. As a corollary we obtain 
that the characteristic classes associated with Cartan connections in 
section \nmb!{2} are the classical characteristic classes of the 
principal bundle $P[H]\to M$.

Section \nmb!{4} deals with a flat Cartan connection of type $\h/\h$ 
on a manifold $P$. We give a simple conceptual proof of the result 
that any flat generalized Cartan connection on a simply connected 
manifold $P$, i\.e\. an $\h$-valued 1-form $\ka$ on $P$ which 
satifies the (left) Maurer-Cartan equation, is the left logalithmic 
derivative of a mapping $\ph:P\to H$ into a Lie group $H$ 
corresponding to $\h$; so `$\ka=\ph\i.d\ph$'. Moreover the mapping 
$\ph$ is uniquely determined up to a left translation. 

A generalized Cartan-connection $\ka:TP\to \h$ induces a homomorphism 
$$\gather
\ka^*:\La(\h^*)\to\Om(P)\\
f\mapsto f\o(\ka\otimes_\wedge \dots\otimes_\wedge \ka)
\endgather$$
of the complex of exterior forms on the Lie algebra $\h$ into the 
complex of differential forms $P$ and (following \cit!{10}) defines a 
a characteristic class of a flat generalized Cartan connection as the 
image of of cohomology classes of the Lie algebra $\h$ under the 
induced homomorphism of cohomologies. This construction may also 
sometimes be applied for the infinite dimensional case.

In section \nmb!{5} we describe a flat Cartan connection associated 
with a flat $G$-structure $p:P\to M$. It defines a Cartan connection 
on the total space $P^\infty$ of the infinite prolongation 
$p^\infty:P^\infty\to M$, which consists of all infinite jets of 
holonomic sections of $p$.

In the last section \nmb!{6} we review shortly the theory of 
prolongation of $G$-structures in the sense of \cit!{22}. Under some 
conditions we define a canonical Cartan connection of type 
$(V\oplus \g^\infty)/\g$ on the total space of the full prolongation 
of a $G$-structure of first or second order.

\head\totoc\nmb0{2}. Cartan connections and generalized Cartan 
connections \endhead

\subhead\nmb.{2.1}. Cartan connections \endsubhead
Let $\h$ be a finite dimensional Lie algebra and let $\g$ be a 
subalgebra of $\h$. Let $P$ be a smooth manifold with $\dim P= 
\dim \h$. By an \idx{\it $\h$-valued absolute parallelism} on $P$ we 
mean a 1-form $\ka\in\Om^1(P;\h)$ with values in 
$\h$ which is non-degenerate in the sense that $\ka_x:T_xP\to \h$ is 
invertible for all $x\in P$. Thus its inverse induces a linear 
mapping $\ze:\h\to \X(P)$ which is given by $\ze_X(x)=(\ka_x)\i(X)$.
Vector fields of the form $\ze_X$ are called \idx{\it parallel}. In 
general, $\ze$ is not  a Lie algebra homomorphism.

\definition{Definition} In this setting
a \idx{\it Cartan connection} of type $\h/\g$ on the manifold $P$ 
is an $\h$-valued absolute parallelism $\ka:TP\to \h$ such that 
\roster
\item $[\ze_X,\ze_Y]=\ze_{[X,Y]}$ for $X\in \h$ and $Y\in \g$. 
\endroster
So the 
inverse mapping $\ze:\h\to\X(P)$ preserves Lie brackets if one of the 
arguments is in $\g$. In particular, the restriction of $\ze$ to $\g$ 
is a Lie algebra homomorphism, and in particular $P$ is a free 
$\g$-manifold.
\enddefinition

\subhead\nmb.{2.2}. Principal Cartan connections on a principal $G$-bundle 
\endsubhead
Let $p:P\to M$ be a principal bundle with structure group $G$ whose 
Lie algebra is $\g$. 
We shall denote by $r:P\x G\to P$ the principal right action and by 
$\ze:\g\to \X(P)$ the fundamental vector field mapping, a Lie 
algebra homomorphism, which is given by $\ze_X(x)=T_e(r_x).X$. Its 
`inverse' is then defined on the vertical bundle $VP$, it is given by 
$\ka_G:VP\to \g$, $\ka_G(\xi_x)=T_e(r_x)\i(\xi_x)$; we call it the 
\idx{\it vertical parallelism}. 

Let us now assume that $\g$ is a subalgebra of a Lie algebra $\h$ 
with $\dim \h = \dim P$.
A $\h/\g$-Cartan 
connection $\ka:TP \to \h$ on $P$ is called a \idx{\it principal 
Cartan connection} of the principal bundle $p:P\to M$, if the 
following two conditions are satisfied:
\roster
\item $\ka|VP=\ka_G$, i\.e\. $\ka$ is an extension of the natural 
       vertical parallelism.
\item $\ka$ is $G$-equivariant, i\.e\. 
       $\ka\o T(r^g)=\Ad(g\i)\o \ka$ 
       for all $g\in G$. If $G$ is connected this follows from 
       \nmb!{2.1},\therosteritem1.
\endroster
This is the usual concept of Cartan connection as used e\.g\. in 
\cit!{13}, p\. 127.

\remark{Remark} Let $\ka\in\Om^1(P;\h)$ be a $\h/\g$-Cartan 
connection on a manifold $P$. Assume that all parallel vector fields 
$\ze(\g)$ are complete. Then they define a locally free action of a 
connected Lie group $G$ with Lie algebra $\g$. If this action is free 
and if the orbit space $M:=P/G$ is a smooth manifold (this is the case 
if the action is also proper), then $p:P\to M$ is a principal 
$G$-bundle and $\ka$ is a principal Cartan connection on it.
\endremark

\subhead\nmb.{2.3}. Principal Cartan connections and a reduction of 
a principal bundle with a connection \endsubhead
Let $H$ be a Lie group with Lie algebra $\h$, let $p:Q\to M$ be a 
principal $H$-bundle, and let $\om:TQ\to \h$ be a principal 
connection form on $Q$. Let us denote by $\Cal H=\ker\om$ the 
horizontal distribution of the connection $\om$. Then we have
$$T_qQ = V_qQ \oplus \Cal H_q,\tag1$$  where $VQ=\ker T(p)\subset TQ$ 
is the vertical subbundle.
We assume now that $G$ is a Lie subgroup of $H$ and that the 
principal bundle $Q\to M$ admits a reduction of the structure group 
to a principal $G$-bundle $p=p|P:P\to M$. So the embedding $P\to Q$ 
is a principal bundle homomorphism over the group homomorphism
$G\to H$:
$$\CD
P @>>> Q \\
@V{p|P}VV @VVpV\\
M @=   M
\endCD$$
Note that for the vertical bundles we have $T_uP\cap V_uQ = V_uP$, 
but the intersection $T_uP\cap \Cal H_u$ may be arbitrary. We have 
the following characterization of the two extremal cases when this 
last intersection is maximal or minimal.

\proclaim{Proposition} \thetag A In the situation above the following 
conditions are equivalent:
\roster
\item For any $u\in P$ the horizontal subspace $\Cal H_u=\ker\om_u$ 
       is contained in $T_uP$, and thus $T_uP = V_uP\oplus \Cal H_u$.
\item The connection $\om$ on $Q$ is induced from a principal 
       connection on $P\to M$ on the associated bundle $Q=P\x_GH$, 
       where $G$ acts on $H$ by conjugation.
\item The holonomy group of the connection $\om$ is contained in $G$.
\endroster
\thetag B The restriction $\om|P= \text{incl}^*\om$ of $\om$ on $P$ 
is a Cartan connection of the principal bundle $p:P\to M$ if and only 
if $T_uP\cap \Cal H_u = 0$ for each $u\in P$, and if 
$\dim M = dim \h/\g$. \qed
\endproclaim

\subhead\nmb.{2.4}. Cartan connections as $G$-structures with 
connections \endsubhead
We establish here a bijective correspondence between principal Cartan 
connections and $G$-structures with a connection.

Let $G\subset GL(V)$, $V=\Bbb R^n$, be a linear Lie group. We recall 
that a $G$-structure on an $n$-dimensional manifold $M$ is a 
principal $G$-bundle $p: P\to M$ together with a displacement form 
$\th:TP\to V$, i\.e\. a $V$-valued 1-form which is $G$-equivariant 
and strictly horizontal in the sense that $\ker\th= VP$.

We assume now that $G$ is a reductive Lie subgroup of a Lie group 
$H$ such that 
$$\h= V\oplus \g,\quad [\g, V]\subset V$$
is the reductive decomposition of the Lie algebra of $H$, and that 
the adjoint representation of $G$ in $V$ is faithful. Then we may 
identify $G$ with a subgroup of $GL(V)$.

\proclaim{Proposition}
In this situation let $\ka: TP\to \h=V\oplus\g$ be a Cartan 
connection on the principal $G$-bundle $p:P\to M$, and let 
$\th=\pr_V\o\ka$ and $\om=\pr_\g\o\ka$ be its components in $V$ and 
$\g$, respectively. 

Then $\th$ is a displacement form and $\om$ is a connection form on 
$p:P\to M$, so that $(p:P\to M,\th)$ is a $G$-structure with a 
connection form $\om$.

Conversely, if $(p:P\to M,\th)$ is a $G$-structure with a connection 
$\om$, then $\ka=\th+\om$ is a principal Cartan connection for the 
principal $G$-bundle $p:P\to M$. \qed
\endproclaim

\subhead\nmb.{2.5}. Generalized Cartan connections \endsubhead
For a principal $G$-bundle $\pi:P\to M$ as in \nmb!{2.2}, if 
$\ka\in\Om^1(P;\h)^G$ is a $G$-equivariant extension of 
$\ka_G:VP\to \g$, we call it a
\idx{\it generalized principal $\h/\g$-Cartan connection}. 

More general, let $P$ be a smooth manifold, let $\h$ be a Lie algebra 
with $\dim \h=\dim P$. We then consider a free action of a Lie 
subalgebra $\g$ of $\h$ on $P$, i\.e\. an injective Lie algebra 
homomorphism $\ze:\g\to\X(P)$. 
A \idx{\it generalized $\h/\g$-Cartan connection} $\ka$ on the 
$\g$-manifold $P$ is then a $\g$-equivariant $\h$-valued one form 
$$\ka\in\Om^1(P;\h)^{\g} :=\{\ph\in\Om^1(P;\h):
     \L_{\ze_X}\ka = \ad(X)\o \ka \text{ for all }X\in \g\}$$ which 
reproduces the generators of the $\g$-fundamental vector fields on 
$P$: for all $X\in\g$ we have $\ka(\ze_X(x))=X$.

\subhead\nmb.{2.6}. Generalized Cartan connections and principal 
connection forms \endsubhead
Let $P$ be a smooth manifold with a free action of a Lie algebra 
$\g$. In \cit!{5} we define the notion of a principal connection on 
$P$ as follows: A principal connection form on $P$ is a $\g$-valued 
$\g$-equivariant 1-form $\om\in\Om(P;\g)^\g$ which reproduces the 
generators of the fundamental vector fields on $P$, so $\om(\ze_X)=X$ 
for $X\in \g$.

As a generalization of proposition \nmb!{2.3} we establish now 
relations between generalized Cartan connections and principal 
connection forms.

\proclaim{Proposition}
Let $\g$ be a reductive subalgebra of a Lie algebra $\h$ with 
reductive decomposition 
$$ \h= V\oplus \g,\quad [\g,V]\subset V.$$
Let $\ka:TP\to \h$ be a generalized Cartan connection on a 
$\g$-manifold $P$ with a free action of the Lie algebra $\g$. 

Then the $\g$-component $\om = \pr_\g\o \ka$ is a principal 
connection form on the $\g$-manifold $P$. In particular, $\ka$ 
defines a $\g$-invariant horizontal distribution $\Cal H:= 
\ka\i(V)\subset TP$ which is a complementary subbundle of the 
`vertical' distribution $\ze_\g(P)\subset P$ spanned by the 
$\g$-action, and which is $\g$-invariant: 
$$[\ze_\g,\Ga(\Cal H)]\subset \Ga(\Cal H),$$
where $\Ga(\Cal H)\subset \X(P)$ is the space of section of the 
bundle $\Cal H$.
\endproclaim

\subhead Remark \endsubhead
It is a natural idea to consider the $V$-component $\th=\pr_V\o\ka$ 
of $\ka$ as some analogon of the notion of displacement form. Clearly 
$\th$ is $\g$-equivariant  and horizontal: 
$\ker\th\supset \ze_\g(P)$. But it will be strictly horizontal ($ 
\ker\th=\ze_\g(P)$) if and only if $\ka$ is a Cartan connection. 
In general we only have $\ker\th = \ze_\g(P) \oplus \Cal K$, where 
$\Cal K = \ker(\ka|\Cal H)$ is a $\g$-invariant distribution, 
possibly of non-constant rank. 

\subhead\nmb.{2.7}. Curvature and Bianchi identity \endsubhead
For a generalized Cartan connection $\ka\in \Om^1(P;\h)^{\g}$ we 
define the \idx{\it curvature} $K$ by
$K= d\ka + \frac12[\ka,\ka]^\wedge$,
where we used the graded Lie bracket on $\Om(P;\h)$ given in 
\cit!{5},~4.1 . From the graded Jacobi identity in $\Om(P;\h)$ we get 
then easily the \idx{\it Bianchi identity}
$$dK + [\ka,K]^\wedge = 0.$$
Then $K$ is \idx{\it horizontal}, i\.e\. kills all $\ze_X$ for 
$X\in\g$, and is $\g$-equivariant, $K\in\Om^2_{\text{hor}}(P;\h)^\g$. 
If $\ka$ is a generalized principal Cartan connection on a principal 
$G$-bundle, then $K$ is even $G$-equivariant, 
$K\in\Om^2_{\text{hor}}(P;\h)^G$.

If $\ka$ is a Cartan connection then an easy computation shows that 
$$\ze^\ka(K(\ze^\ka_X,\ze^\ka_Y))
     =[\ze^\ka_X,\ze^\ka_Y]_{\X(P)}-\ze^\ka([X,Y]_\h).$$

\subhead\nmb.{2.8}. Covariant exterior derivative \endsubhead
For a generalized $\h/\g$-Cartan connection $\ka\in\Om^1(P;\h)^\g$ and any 
representation $\rh:\h\to GL(W)$ we define the 
\idx{\it covariant exterior derivative} 
$$\gather
d_\ka:\Om^p_{\text{hor}}(P;W)^\g \to \Om^{p+1}_{\text{hor}}(P;W)^\g \\
d_\ka\Ps = d\Ps + \rh^\wedge (\ka)\Ps.
\endgather$$
For a principal Cartan connection on a principal $G$-bundle we even
have
$$d_\ka(\Om^p_{\text{hor}}(P;W)^G)\subset\Om^{p+1}_{\text{hor}}(P;W)^G.$$

\subhead\nmb.{2.9}. Chern-Weil forms \endsubhead
If $f\in L^k(\h) := (\bigotimes^k\h^*)$ is a $k$-linear function 
on $\h$ and if $\ps_i\in\Om^{p_i}(P;\h)$ 
we can construct the following differential forms 
$$\gather
\ps_1\otimes_\wedge \dots \otimes_\wedge \ps_k 
	\in \Om^{p_1+\dots+p_k}(P;\h\otimes \dots \otimes \h),\\
f^{\ps_1,\dots,\ps_k}:= f\o (\ps_1\otimes_\wedge \dots \otimes_\wedge \ps_k ) 
	\in \Om^{p_1+\dots+p_k}(P).
\endgather$$
The exterior derivative of the latter one is clearly given by
$$\multline
d( f\o(\ps_1\otimes_\wedge\dots\otimes_\wedge\ps_k)) = 
	f\o d(\ps_1\otimes_\wedge\dots\otimes_\wedge\ps_k) =\\
= f\o \left( \tsize\sum_{i=1}^k(-1)^{p_1+\dots+p_{i-1}} 
  	\ps_1\otimes_\wedge\dots\otimes_\wedge d\ps_i\otimes_\wedge
	\dots\otimes_\wedge\ps_k \right).
\endmultline$$
Note that the form $f^{\ps_1,\dots,\ps_k}$ is $\g$-invariant and 
horizontal if all $\ps_i\in \Om^{p_i}_{\text{hor}}(P;\h)^\g$ and 
$f\in L^k(\h)^\g$ is invariant under the adjoint action of $\g$ on $\h$. 
It is then the pullback of a form on $M$. For a principal Cartan 
connection one may replace $\g$ by $G$.

\proclaim{\nmb.{2.10}. Lemma} Let $\ka$ be a generalized $\h/\g$-Cartan 
connection on $P$.
Let $f\in L^k(\h)^\h$ be $\h$-invariant under the adjoint action then 
the differential form $f^K:=f^{K,\dots,K}$ is closed in 
$\Om_{\text{hor}}^{2k}(M)^\g$.
\endproclaim

\demo{Proof}
The same computation as in the proof of \cit!{5},~7.4  with $\om$ and 
$\Om$ replaced by $\ka$ and $K$.
\qed\enddemo

\proclaim{\nmb.{2.11}. Proposition} Let $\ka_0$ and $\ka_1$ be two 
generalized $\h/\g$-Cartan connections on $P$ with curvature forms $K_0, 
K_1\in \Om^2(P;\h)^\g$, and let $f\in L^k(\h)^\h$. Then the cohomology 
classes of the two closed forms $f^{K_0}$ and $f^{K_1}$ in 
$H^{2k}(\Om^*_{\text{hor}}(P)^{\g})$ agree.

If $P\to M$ is a principal $G$-bundle and if $\ka_1$ and $\ka_2$ are 
principal generalized Cartan connections on it, then the cohomology 
classes of the two closed forms $f^{K_0}$ an $f^{K_1}$ agree in 
$H^{2k}(M)$.
\endproclaim

\demo{Proof}
Literally the same proof as for \cit!{5},~7.5 applies, with $\om$ and 
$\Om$ replaced by $\ka$ and $K$.
\qed\enddemo

\head\totoc\nmb0{3}. The relation between principal Cartan connections and 
principal connections \endhead

In this section we follow the notation and concepts of \cit!{14}, 
chapter III, which we also explain here.

\subhead\nmb.{3.1}. Extension of the structure group \endsubhead
Given a principal bundle $\pi:P\to M$ with structure group $G$ and 
$G\subset H$ we consider the left action of $G$ on $H$ (by left 
translation) and the associated bundle $\pi:P[H]=P\x_GH\to M$. Recall 
from \cit!{14},~10.7 the $G$-bundle projection $q:P\x H\to P[H]=P\x_GH$.
Since $q(u.g,h)=q(u,gh)$ we get 
$Tq(Tr(X_u,Z_g),Y_h)=Tq(X_u,T\la_g.Y_h+T\rh_h.Z_g)$.
This is then a principal $H$-bundle with principal $H$-action 
$\tilde r: P[H]\x H\to P[H]$ given by $\tilde r(q(u,h),h')=q(u,hh')$. 
Since $G\subset H$ is $G$-equivariant we get a homomorphism of 
principal bundles over $G\subset H$
$$\newCD
P @(2,0) @()\l{\pi}@(1,-1) & & P[H] @()\l{\pi}@(-1,-1) \\
& M &
\endnewCD$$

\proclaim{\nmb.{3.2}. Lemma} In the situation of \nmb!{3.1} the 
generalized Cartan connections in the space $\Om^1(P;\h)^G$ correspond 
canonically and bijectively to the $H$-principal connections in 
$\Om^1(P[H];\h)^H$.
\endproclaim

\demo{Proof}  For $Y\in\h$ the fundamental vector field 
$\ze^{P[H]}_Y$ on $P{H}$ is given by 
$$\ze^{P[H]}_Y(q(u,h))=T(\tilde r)(Tq(0_u,0_h),Y) = 
T_{(u,h)}q(0_u,T\la_h.Y).$$
For a generalized Cartan connection $\ka\in\Om^1(P;\h)^G$ we define 
for $X_u\in T_uP$ and $Y\in \h$:
$$\gather
(q^\flat\ka):TP[H]\to \h\\
(q^\flat\ka)(T_{(u,h)}q(X_u,T_e\la_h.Y)):= Y + 
     \Ad(h\i)\ka_u(X_u)\tag1
\endgather$$
We claim that $(q^\flat\ka)\in\Om^1(P[H];\h)^H$ is well defined and is 
a principal connection.

It is well defined: for $g\in G$ we have
$$\align
(q^\flat\ka)(T_{(u,h)}q(Tr^g.X_u,T\la_{g\i}T_e\la_h.Y)) &= 
     Y + \Ad((g\i h)\i)\ka_u(Tr^g.X_u) \\
&=   Y + \Ad(h\i)\Ad(g)\Ad(g\i)\ka_u(X_u) \\
&= (q^\flat\ka)(T_{(u,h)}q(X_u,T_e\la_h.Y))
\endalign$$
Moreover $T_{(u,h)}q(X_u,Y_h)=0$ if and only if 
$(X_u,Y_h)=(\ze^P_X(u),-T\rh_h.X)$ for some $X\in\g$, but then
$$\align
(q^\flat\ka)(\ze^P_X(u),-T\rh_h.X) &= 
     (q^\flat\ka)(\ze^P_X(u),-T_e\la_h.\Ad(h\i)X) \\
&= -\Ad(h\i)X + \Ad(h\i)\ka_u(\ze^P_X(u)) = 0.
\endalign$$
We check that it is $H$-equivariant:
$$\align
(q^\flat\ka)(T(\tilde r^k).T_{(u,h)}q.(X_u,T_e\la_h.Y)) 
     &= (q^\flat\ka)(T_{(u,hk)}q.(X_u,T\rh_k.T_e\la_h.Y)) \\
&= (q^\flat\ka)(T_{(u,hk)}q.(X_u,T_e\la_{hk}.T\rh_k.T\la_{k\i}.Y)) \\
&= (q^\flat\ka)(T_{(u,hk)}q.(X_u,T_e\la_{hk}.\Ad(k\i)Y)) \\
&= \Ad(k\i)Y + \Ad(k\i.h\i)\ka_u(X_u) \\
&= \Ad(k\i)(q^\flat\ka)(T_{(u,h)}q.(X_u,T_e\la_h.Y)).
\endalign$$ 
Next we check that it reproduces the infinitesimal generators of 
fundamental vector fields:
$$ (q^\flat\ka)(\ze^{P[H]}_Y(q(u,h)) 
     = (q^\flat\ka)(T_{(u,h)}q(0_u,T\la_h.Y))=Y$$

Now let $\om\in\Om^1(P[H];\h)^H$ be a principal connection form.
Then the pull back $(q^\flat)\i\om$ of $\om$ to the $G$-subbundle 
$P\subset p[H]$ is in $\Om^1(P;\h)^G$ and clearly reproduces the 
infinitesimal generators of $G$-fundamental vector fields, so it is a 
generalized Cartan connection.
Explicitely we have $((q^\flat)\i\om)(X_u):=\om(T_{(u,e)}q(X_u,0_e))$
and with this formula it is easy to check that the two construction
are inverse to each other.
\qed\enddemo

\proclaim{\nmb.{3.3}. Theorem} Let $\pi:P\to M$ be a principal 
bundle with structure group $G$, let $H$ be a Lie group containing 
$G$ and let $\rh:H\to GL(W)$ be a finite dimensional linear 
representation of $H$.

Then there is a canonical linear isomorphism
$$q^\flat:\Om^p_{\text{hor}}(P;W)^G\to\Om^p_{\text{hor}}(P[H];W)^H$$
which intertwines the covariant exterior derivatives of any 
generalized Cartan connection $\ka$ on $P$ with values in $\h$ and of 
its canonically associated principal connection $q^\flat\ka$ on 
$P[H]$:
$$d_{q^\flat\ka}\o q^\flat = q^\flat\o d_\ka:
     \Om^p_{\text{hor}}(P;W)^G \to \Om^{p+1}_{\text{hor}}(P[H];W)^G$$
If $K\in \Om^2_{\text{hor}}(P;\h)^G$ is the curvature of a generalized 
Cartan connection $\ka$ in the sense of \nmb!{2.7} then 
$q^\flat K\in\Om^2_{\text{hor}}(P[H];\h)^H$ is the principal curvature 
of the principal connection $q^\flat\ka$ on $P[H]$.
\endproclaim

\demo{Proof}
For $\Ps\in\Om^p_{\text{hor}}(P;\h)^G$ we define 
$q^\flat\Ps\in\Om^p_{\text{hor}}(P[H];\h)^H$ by
$$\align 
(q^\flat\Ps)_{q(u,h)}&(Tq(\xi^1_u,T_e\la_h.Y^1),
     Tq(\xi^2_u,T_e\la_h.Y^2),\dotsc)=\tag1\\
&=\Ad(h\i)\Ps_u(\xi^1_u,\xi^2_u,\dotsc).
\endalign$$
This is well defined and horizontal, since a vector 
$Tq(\xi_u,T\la_h.Y)$ is vertical in $P[H]$ if and only if it is of the 
form $Tq(\ze^P_X(u),T_e\la_h(Z-\Ad(h\i)X))$ for some $Z\in\h$ and 
$X\in\g$, and the right hand side vanishes if one such vector appears 
in the left hand side. Note that $q^\flat\Ps$ is well defined only if 
$\Ps$ is horizontal. It is easily seen that $q^\flat\Ps$ is 
$H$-equivariant.

If $\Ph\in\Om^p_{\text{hor}}(P[H];\h)^H$ then the pull back of $\Ph$ 
to the subbundle $P$ gives a form 
$(q^\flat)\i\Ph\in\Om^p_{\text{hor}}(P;\h)^G$.
We have the explicit formula 
$((q^\flat)\i\Ph)(\xi^1_u, \xi^2_u,\dotsc) 
     = \Ph(Tq(\xi_u^1,0_e),Tq(\xi^2_u,0_e),\dotsc)$,
and using this it is easy to show that the two constructions are 
inverse to each other:
$$\align
((q^\flat)(q^\flat)\i&\Ph)_{q(u,h)}(Tq(\xi^1_u,T\la_h.Y^1),\dotsc)
     = \Ad(h\i)((q^\flat)\i\Ph)_u(\xi^1_u,\xi^2_u,\dotsc)\\
&= \Ad(h\i)\Ph(T_{(u,e)}q(\xi^1_u,0_e),\dotsc)
     = \Ph(T(\tilde r^h).T_{(u,e)}q(\xi^1_u,0_e),\dotsc)\\
&= \Ph(T_{(u,h)}q(\xi^1_u,T\la_h.0_e),\dotsc)
     = \Ph(T_{(u,h)}q(\xi^1_u,T\la_h.Y),\dotsc),
\endalign$$
since $\Ph$ is horizontal, and
$$((q^\flat)\i(q^\flat)\Ps)_u(\xi^1_u,\xi^2_u,\dotsc)
     =((q^\flat)\i\Ps)_{q(u,e)}(Tq(\xi^1_u,0_e),\dotsc)
     =\Ps_u(\xi^1_u,\dotsc).$$

{\bf Claim 1:}
$d_{q^\flat\ka}\o q^\flat = q^\flat\o d_\ka:
     \Om^p_{\text{hor}}(P;W)^G \to \Om^{p+1}_{\text{hor}}(P[H];W)^G$ 
holds for a generalized Cartan connection $\ka$ on $P$. 
Here $d_{q^\flat\ka}$ is given by $d_{q^\flat\ka}\Ph=\ch^*d\Ph$ for any 
form $\Ph\in\Om(P[H],V)$ with values in a vector space $V$, where 
$\ch$ is the horizontal projection induced by ${q^\flat\ka}$. In 
\cit!{4},~1.4 it is proved that for $\Ph\in\Om_{\text{hor}}(P[H];W)^H$ 
the formula $d_{q^\flat\ka}\Ph=d\Ph+[q^\flat\ka,\Ph]^\wedge$ holds.
On the other hand we have $d_\ka\Ps = d\Ps + \rh^\wedge (\ka)\Ps$ for 
$\Ps\in\Om_{\text{hor}}(P;W)^G$ by definition \nmb!{2.8}.

To compute $d(q^\flat\Ps)$ we need vector fields. So let 
$\xi_i\in\X(P)^G$ be $G$-equivariant vector fields on $P$, and for 
$Y_i\in\h$ let $L_{Y_i}$ denote the left invariant vector field on $H$, 
$L_{Y_i}(h)=T\la_h.Y_i$. Then the vector field $\xi_i\x L_{Y_i}$ is 
$G$-equivariant and factors thus to a vector field on the associated 
bundle as indicated in the following diagram:
$$\cgaps{.5;1.5;.5}\newCD
& P\x H @()\L{\xi_i\x L_{Y_i}}@(1,0) @()\L{q}@(0,-1) & 
TP\x TH @()\l{Tq}@(0,-1) & \\
P[H] @()\a=@(1,0) & P\x_GH @()\L{\widetilde{\xi_i\x L_{Y_i}}}@(1,0) &
TP\x_{TG}TH @()\a=@(1,0) & T(P[H])
\endnewCD$$
So the vector fields $\xi_i\x L_{Y_i}$ on $P\x H$ and 
$\widetilde{\xi_i\x L_{Y_i}}$ on $P[H]$
are $q$-related and thus we have
$$[\widetilde{\xi_i\x L_{Y_i}},\widetilde{\xi_j\x L_{Y_j}}] = 
\widetilde{[\xi_i,\xi_j]\x L_{[Y_i,Y_j]}}\tag2$$
Now we compute 
$$\align
d(q^\flat\Ps)&\left(\widetilde{\xi_0\x L_{Y_0}},\dots,
     \widetilde{\xi_p\x L_{Y_p}}\right) = \\
&= \sum_{i=0}^p(-1)^i (\widetilde{\xi_i\x L_{Y_i}})\left(q^\flat\Ps\left(
     \widetilde{\xi_0\x L_{Y_0}},\dots,\widehat{\quad_i},\dots,
     \widetilde{\xi_p\x L_{Y_p}}  \right)\right)\\
&\quad + \sum_{i<j}(-1)^{i+j}(q^\flat\Ps)\left(\left
     [\widetilde{\xi_i\x L_{Y_i}},\widetilde{\xi_i\x L_{Y_i}}\right],
     \widetilde{\xi_0\x L_{Y_0}},\dots,\widehat{\quad_i},
     \dots,\widehat{\quad_j},\dotsc\right).
\endalign$$
Since we have 
$$\align
(q^\flat\Ps)_{q(u,h)}&\left(\widetilde{\xi_1\x L_{Y_1}},\dots,
     \widetilde{\xi_p\x L_{Y_p}}\right) =\tag3\\
&= (q^\flat\Ps)_{q(u,h)}(Tq(\xi_1(u),T_e\la_h.Y_1),\dots,
     Tq(\xi_p(u),T_e\la_h.Y_p))\\
&= \Ad(h\i).\Ps_u(\xi_1(u),\dots,\xi_p(u))\in \h
\endalign$$
we get
$$\align
(\widetilde{\xi_0\x L_{Y_0}})&\left(q^\flat\Ph\left(
     \widetilde{\xi_1\x L_{Y_1}},\dots,\widetilde{\xi_p\x L_{Y_p}}
     \right)\right)(q(u,h))=\\
&=(T_h(\Ad\o\operatorname{Inv}).T_e\la_h.Y_0).\Ps_u(\xi_1,\dots,\xi_p) +
     \Ad(h\i)(\xi_0\Ps(\xi_1,\dots,\xi_p))\\
&=-[Y_0,\Ad(h\i).\Ps_u(\xi_1,\dots,\xi_p)] + 
     \Ad(h\i)(\xi_0\Ps(\xi_1,\dots,\xi_p)).
\endalign$$
Inserting we get
$$\align
d(q^\flat\Ps)&\left(\widetilde{\xi_0\x L_{Y_0}},\dots,
     \widetilde{\xi_p\x L_{Y_p}}\right)(q(u,h)) = \\
&= - \sum_{i=0}^p(-1)^i[Y_i,\Ad(h\i).\Ps_u(\xi_0,\dots, 
     \widehat{\xi_i},\dots\xi_p)] +\\
&\quad + \Ad(h\i).(d\Ps)_u(\xi_0,\dots,\xi_p).
\endalign$$
Next we compute 
$$\align
[q^\flat\ka,&q^\flat\Ps]^\wedge \left(\widetilde{\xi_0\x L_{Y_0}},\dots,
     \widetilde{\xi_p\x L_{Y_p}}\right)(q(u,h)) = \\
&= \sum_{i=0}^p(-1)^i \left[(q^\flat\ka)_{q(u,h)}
     \left(\widetilde{\xi_i\x L_{Y_i}}\right),
     (q^\flat\Ps)_{(q(u,h))}\left(\widetilde{\xi_0\x L_{Y_0}},\dots,
     \widehat{\quad_i},\dotsc\right)\right]\\
&= \sum_{i=0}^p(-1)^i [Y_i+\Ad(h\i)\ka_u(\xi_i),
     \Ad(h\i)\Ps_u(\xi_0,\dots,\widehat{\xi_i},\dots,\xi_p)]_\h\\
&= \sum_{i=0}^p(-1)^i[Y_i,\Ad(h\i).\Ps_u(\xi_0,\dots, 
     \widehat{\xi_i},\dots\xi_p)] +\\
&\quad + \Ad(h\i).[\ka,\Ps]^\wedge (\xi_0,\dots,\xi_p)(u).
\endalign$$
On the other hand we have 
$$\align
(q^\flat\o d_\ka\Ps)_{(q(u,h))}\left(\widetilde{\xi_0\x L_{Y_0}},\dots,
     \widetilde{\xi_p\x L_{Y_p}}\right) 
     &=\Ad(h\i).(d_\ka\Ps)_u(\xi_0,\dots,\xi_p)\\
&=\Ad(h\i).(d\Ps+[\ka,\Ps])_u(\xi_0,\dots,\xi_p),
\endalign$$
so the claim follows by comparing the last three expressions.

{\bf Claim 2:} 
$q^\flat K = d(q^\flat\ka)+\frac12[q^\flat\ka,q^\flat\ka]^\wedge$
for a generalized Cartan connection $\ka$ on $P$ with curvature 
$K=d\ka+\frac12[\ka,\ka]^\wedge$. 

We have $K\in\Om^2_{\text{hor}}(P;\h)^G$ but $\ka$ is not horizontal, 
so we must redo parts of the above computations. We use the same 
vector fields as in the proof of Claim 1. Since by 
\nmb!{3.2},\thetag1 we have
$$\align
(q^\flat\ka)_{q(u,h)}\left(\widetilde{\xi_1\x L_{Y_1}}\right) 
&= (q^\flat\ka)_{q(u,h)}(Tq(\xi_1(u),T_e\la_h.Y_1))\tag4\\
&= Y_1 + \Ad(h\i).\ka_u(\xi_1)
\endalign$$
we get again the same formula as for $\Ps$
$$\align
\widetilde{\xi_0\x L_{Y_0}}
     &(q^\flat\ka)_{q(u,h)}\left(\widetilde{\xi_1\x L_{Y_1}}\right)\\
&= -[Y_0,\Ad(h\i)\ka_u(\xi_1)] + \Ad(h\i)(\xi_0\ka(\xi_1))(u).
\endalign$$
This leads to 
$$\align
d(q^\flat\ka)&\left(\widetilde{\xi_0\x L_{Y_0}},
     \widetilde{\xi_1\x L_{Y_1}}\right) = \\
&= -[Y_0,\Ad(h\i)\ka_u(\xi_1)] + [Y_1,\Ad(h\i)\ka_u(\xi_0)] - [Y_0,Y_1]\\
&\quad + \Ad(h\i)(d\ka_u(\xi_0,\xi_1)).
\endalign$$
Again from \thetag4 we get
$$\align
\tfrac12[q^\flat\ka,q^\flat\ka]^\wedge &\left(\widetilde{\xi_0\x L_{Y_0}},
     \widetilde{\xi_1\x L_{Y_1}}\right)(q(u,h)) = \\
&= \tfrac12[Y_0+\Ad(h\i)\ka_u(\xi_0),Y_1+\Ad(h\i)\ka_u(\xi_1)] \\
&\quad - \tfrac12[Y_1+\Ad(h\i)\ka_u(\xi_1),Y_0+\Ad(h\i)\ka_u(\xi_0)]\\
&= [Y_0,\Ad(h\i)\ka_u(\xi_1)] - [Y_1,\Ad(h\i)\ka_u(\xi_0)] + [Y_0,Y_1]\\
&\quad + \tfrac12\Ad(h\i).[\ka,\ka]^\wedge (\xi_0,\xi_1)(u)
\endalign$$
from which now the result follows.
\qed\enddemo

\proclaim{\nmb.{3.4}. Corollary} The characteristic class for an 
invariant $f\in L^k(\h)^H$ constructed in proposition \nmb!{2.11} with 
the help of generalized Cartan connections on $P$ is exactly the 
characteristic class of the principal bundle $P[H]$ associated to 
$f$. Since $P[H]$ admits a reduction of the structure group to 
$G$, this class is a characteristic class of $P$, associated to 
$f|\g\in L^k(\g)^G$. If $f|\g=0$ then the form $f^K$ of proposition 
\nmb!{2.11} is exact.
\endproclaim

\demo{Proof}
This follows from well known properties of characteristic classes of 
principal bundles.
\qed\enddemo

\head\totoc\nmb0{4}. Flat Cartan connections \endhead

\subhead\nmb.{4.1}. Flat Cartan connections \endsubhead
Let $P$ be a smooth manifold. A Cartan connection $\ka:TP\to\h$ is 
said to be flat if its curvature $K=d\ka+\frac12[\ka,\ka]^\wedge $ 
(see \nmb!{2.7}) vanishes. In this case the subalgebra $\g\subset \h$ 
does not play any role. The inverse mapping 
$\ze:\h\to \X(P)$, given by $\ze_X(x)=(\ka_x)\i(X)$ is then a 
homomorphism of Lie algebras, and it defines a free transitive 
action of the Lie algebra $\h$ on the manifold $P$ in the sense of 
\cit!{5},~2.1. The inverse statement is also valid, see 
\cit!{5},~5.1.

A flat generalized Cartan connection is then a form $\ka:TP\to \h$ 
which satisfies the Maurer-Cartan equation 
$d\ka+\frac12[\ka,\ka]^\wedge$ (without the assumtion that it is 
non-degenerate).

\subhead\nmb.{4.2} \endsubhead
Let $H$ be a connected Lie group with 
Lie algebra $\h$, multiplication $\mu:H\x H\to H$, and for $g\in H$ 
let $\mu_g, \mu^g:H\to H$ denote the left and right translation, 
$\mu(g,h)=g.h=\mu_g(h)=\mu^h(g)$.  
For a smooth mapping $\ph:P\to H$ let us use the left trivialization of 
$TH$ and consider the 
\idx{\it left logarithmic derivative} $\de^l \ph\in\Om^1(P;\h)$, given 
by $\de^l \ph_x := T(\mu_{\ph(x)\i})\o T_x\ph:T_xP\to T_{\ph(x)}H\to \h$. 
Similarly we consider the 
\idx{\it right logarithmic derivative} $\de^r \ph\in\Om^1(P;\h)$ which 
is given by 
$\de^r \ph_x:=T(\mu^{\ph(x)\i})\o T_x\ph:T_xP\to T_{\ph(x)}H\to \h$.
The following result can be found in \cit!{17}, \cit!{18}, 
\cit!{19}, or in \cit!{11} (proved with moving 
frames); see also \cit!{5},~5.2. 
We include a simple conceptual proof and we consider all variants. 

\proclaim{Proposition} For a smooth mapping $\ph:P\to H$ the left
logarithmic derivative $\de^l \ph\in\Om^1(P;\h)$ satisfies the (right)
Maurer-Cartan equation $d\de^l \ph+\frac12[\de^l \ph,\de^l \ph]^\wedge =0$.

If conversely a 1-form $\ka\in\Om^1(P;\h)$ satisfies 
$d\ka+\frac12[\ka,\ka]^\wedge=0$ then for each simply connected subset 
$U\subset P$ there exists a smooth function $\ph:U\to H$ with 
$\de^l \ph=\ka|U$, and $\ph$ is uniquely detemined up to a right translation 
in $H$.

For a smooth mapping $\ph:P\to H$ the right
logarithmic derivative $\de^r \ph\in\Om^1(P;\h)$ satisfies the (left)
Maurer-Cartan equation $d\de^r \ph-\frac12[\de^r \ph,\de^r \ph]^\wedge =0$.

If a 1-form $\ka\in\Om^1(P;\h)$ satisfies 
$d\ka-\frac12[\ka,\ka]^\wedge=0$ then for each simply connected subset 
$U\subset P$ there exists a smooth function $\ph:U\to H$ with 
$\de^r \ph=\ka|U$, and $\ph$ is uniquely determined up to a left translation 
in $H$.
\endproclaim

\demo{Proof} Let us treat first the right logarithmic derivative 
since it leads to a principal connection for a bundle with right 
principal action.
We consider the trivial principal bundle 
$\operatorname{pr_1}:P\x H\to P$ with right principal action. Then 
the submanifolds $\{(x,\ph(x).g):x\in P\}$ for $g\in H$ form a 
foliation of $P\x G$ whose tangent distribution is transversal to the 
vertical bundle $P\x TH \subset T(P\x H)$ and is invariant under the 
principal right $H$-action. So it is the horizontal distribution of a 
principal connection on $P\x H\to H$. For a tangent vector 
$(X_x,Y_g)\in T_xP\x T_gH$ the horizontal part is the right translate 
to the foot point $(x,g)$ of $(X_x,T_x\ph.X_x)$, so the decomposition in 
horizontal and vertical parts according to this distribution is
$$
(X_x,Y_g) = (X_x, T(\mu^g).T(\mu^{\ph(x)\i}).T_x\ph.X_x)
+(0_x,Y_g-T(\mu^g).T(\mu^{\ph(x)\i}).T_x\ph.X_x).
$$
Since the fundamental vector fields for the right action on $H$ are 
the left invariant vector fields, the corresponding connection form 
is given by
$$\align
\om^r(X_x,Y_g) &= 
T(\mu_{g\i}).(Y_g-T(\mu^g).T(\mu^{\ph(x)\i}).T_x\ph.X_x),\\
\om^r_{(x,g)} &= T(\mu_{g\i}) - \Ad(g\i).\de^r \ph_x,\\
\om^r &= \ka_H^l - (\Ad\o\operatorname{Inv}).\de^r \ph,\tag1
\endalign$$
where $\ka_H^l:TH\to\h$ is the left Maurer-Cartan form on $H$ (the 
left trivialization), given by $(\ka^l_H)_g=T(\mu_{g\i})$. Note that 
$\ka^l_H$ is the principal connection form for the (unique) principal 
connection $p:H\to \text{point}$ with right principal action, which 
is flat so that the right (from right action) Maurer-Cartan 
equation holds in the form
$$d\ka^l+\tfrac12[\ka^l,\ka^l]^\wedge=0.\tag2$$

The principal connection $\om^r$ is flat since we got it via the horizontal 
leaves, so the principal curvature form vanishes:
$$\align
0&= d\om^r+\tfrac12[\om^r,\om^r]^\wedge\tag3\\ 
&= d\ka^l_H + \tfrac12[\ka^l_H,\ka^l_H]^\wedge 
     - d(\Ad\o\operatorname{Inv}) \wedge \de^r \ph - 
     (\Ad\o\operatorname{Inv}). d\de^r \ph \\
&\quad - [\ka^l_H,(\Ad\o\operatorname{Inv}). \de^r \ph]^\wedge + 
     \tfrac12[(\Ad\o\operatorname{Inv}).\de^r \ph,
     (\Ad\o\operatorname{Inv}). \de^r \ph]^\wedge\\
&= - (\Ad\o\operatorname{Inv}). (d\de^r \ph - 
     \tfrac12[\de^r \ph,\de^r \ph]^\wedge), 
\endalign$$ 
where we used \thetag2
and since for $X\in\g$ we have:
$$\align
d(\Ad\o\operatorname{Inv}) (T(\mu_g)X) &= \ddt \Ad(\exp(-tX).g\i) 
     = -\ad(X)\Ad(g\i) \\ 
&= -\ad(\ka_H^l(T(\mu_g)X))(\Ad\o\operatorname{Inv})(g),\\
d(\Ad\o\operatorname{Inv}) 
&=-(\ad\o\ka_H^l)(\Ad\o\operatorname{Inv}).\tag4 
\endalign$$
So we have $d\de^r \ph - \tfrac12[\de^r \ph,\de^r \ph]^\wedge$ as asserted.

If conversely we are given a 1-form $\ka^r\in\Om^1(P;\h)$ with 
$d\ka^r - \tfrac12[\ka^r,\ka^r]^\wedge=0$ then we consider the 
1-form $\om^r\in\Om^1(P\x H;\h)$, given by the analogon of \thetag1,
$$
\om^r = \ka_H^l - (\Ad\o\operatorname{Inv}).\ka^r\tag5
$$
Then $\om^r$ is a principal connection form on $P\x H$, since it 
reproduces the generators in $\h$ of the fundamental vector fields 
for the principal right action, i\.e\. the left invariant vector 
fields, and $\om^r$ is $H$-equivariant:
$$\align
((\mu^g)^*\om^r)_h &= \om^r_{hg}\o (Id\x T(\mu^g)) = 
     T(\mu_{g\i.h\i}).T(\mu^g) - \Ad(g\i.h\i).\ka^r\\
&= \Ad(g\i).\om^r_h.
\endalign$$
The computation in \thetag3 for $\ka^r$ instead of $\de^r \ph$ shows 
that this connection is flat. So the horizontal bundle is integrable, 
and $\operatorname{pr_1}:P\x H\to P$, restricted to each horizontal 
leaf, is a covering. Thus it may be inverted over each simply 
connected subset $U\subset P$, and the inverse $(Id,\ph):U\to P\x H$ is 
unique up to the choice of the branch of the covering, and the choice 
of the leaf, i\.e\. $\ph$ is unique up to a right translation by an 
element of $H$. The beginning of this proof then shows that 
$\de^r \ph=\ka^r|U$. 

For the left logarithmic derivative $\de^l \ph$ the proof is similar, 
and we discuss only the essential deviations. First note that on the 
trivial principal bundle $\operatorname{pr_1}:P\x H\to P$ with left 
principal action of $H$ the fundamental vector fields are the right 
invariant vector fields on $H$, and that for a principal connection 
form $\om^l$ the curvature form is given by 
$d\om^l-\frac12[\om^l,\om^l]^\wedge$. Look at the proof of 
\cit!{14},~11.2 to see this. The connection form is then given by
$$
\om^l = \ka_H^r - \Ad.\de^l \ph,\tag1'
$$
where the right Maurer-Cartan form 
$(\ka_H^r)_g=T(\mu^{g\i}):T_gH\to\h$ now satifies the left 
Maurer-Cartan equation
$$d\ka_H^r-\frac12[\ka_H^r,\ka_H^r]^\wedge = 0. \tag2'$$
Flatness of $\om^l$ now leads to the computation
$$\align
0&= d\om^l-\tfrac12[\om^l,\om^l]^\wedge\tag3'\\ 
&= d\ka^r_H - \tfrac12[\ka^r_H,\ka^r_H]^\wedge 
     - d\Ad \wedge \de^l \ph - \Ad. d\de^l \ph \\
&\quad + [\ka^r_H,\Ad. \de^l \ph]^\wedge - 
     \tfrac12[\Ad.\de^l \ph,\Ad. \de^l \ph]^\wedge\\
&= - \Ad. (d\de^l \ph + \tfrac12[\de^l \ph,\de^l \ph]^\wedge), 
\endalign$$ 
where we used
$$\align
d\Ad (T(\mu^g)X) &= \ddt \Ad(\exp(tX).g) = \ad(X)\Ad(g) \\ 
&= \ad(\ka_H^r(T(\mu^g)X))\Ad(g),\\
d\Ad &=(\ad\o\ka_H^r)\Ad.\tag4' 
\endalign$$
The rest of the proof is obvious.
\qed\enddemo

\subhead\nmb.{4.3}. Characteristic classes for flat Cartan 
connections \endsubhead
A generalized Cartan connection $\ka:TP\to \h$ on the manifold $P$ 
induces a homomorphism 
$$\gather
\ka^*:\La(\h^*)\to \Om(P),\\
f\mapsto f^\ka = f\o (\ka\otimes_\wedge \dots\otimes_\wedge \ka)
\endgather$$
of the algebra of exterior forms on $\h$ into the algebra of 
differential forms on $P$. Let us assume now that the Cartan 
connection $\ka$ is flat. Then $\ka^*$ commutes with the exterior 
differentials and is a homomorphism of differential complexes, since 
we have by \nmb!{2.9}
$$\align
d(f(\ka,\dots,\ka))&= 
     \sum_{i=1}^k(-1)^{i-1}f(\ka,\dots,d\ka,\dots,\ka)\\ 
&=\sum_{i=1}^k(-1)^{i-1}f(\ka,\dots,-\tfrac12[\ka,\ka]^\wedge ,\dots,\ka)\\ 
&= (df)(\ka,\dots,\ka)\quad k+1\text{ times}.
\endalign$$
Thus we have an associated homomorphism $\ka^*:H^*(\h,\Bbb R) \to H^*(P)$. 
The nontrivial elements of its image are called 
\idx{\it characteristic classes} for the flat Cartan connection 
$\ka$. In \cit!{6} a similar construction is applied even for 
infinite dimensional manifolds.

\subhead\nmb.{4.4} \endsubhead
Let $\ka:TP\to \h$ be a flat generalized Cartan connection on the 
manifold $P$, so $d\ka+\frac12[\ka,\ka]^\wedge = 0$ holds. Let $H$ be 
a connected Lie group with Lie algebra $\h$. Suppose that there 
exists a smooth mapping $\ph:P\to H$ with $\de^l\ph=\ka$. By 
proposition \nmb!{4.2} such a $\ph$ is unique up to a right 
translation in $H$, and it exists if e\.g\. $P$ is simply connected. 
Clearly $\ph$ is a local diffeomorphism if and only if $\ka$ is a 
Cartan connection (non degenerate). Such a mapping $\ph$ is called 
\idx{\it Cartan's developing}, see also \cit!{5},~5.2. It gives a 
convenient way to express the characteristic classes of \nmb!{4.3}, 
by the following easy results. 

\proclaim{Lemma} Let $\ka:TP\to \h$ be a flat generalized Cartan 
connection such that a Cartan's developing $\ph:P\to H$ exists.
Then the following diagram commutes
$$\newCD
\La^k(\h^*) @()\L{\ka^*}@(2,0) @()\l{L}@(1,-1) & & \Om^k(P)\\
 & \Om^k(H) @()\l{\ph^*}@(1,1) &
\endnewCD$$
where $L$ is the extension to left invariant differential forms.
\endproclaim

\demo{Proof}
Plug in the definitions.
\qed\enddemo

\subhead\nmb.{4.5}. Cartan connections with constant curvature
\endsubhead
Let $\ka:TP\to \h$ be a Cartan connection of type $\h/\g$ on a 
manifold $P$. Then its curvature belongs to the space 
$\Om^2(P;\h)\cong \h\otimes \Om^2(P)$. Using the absolute 
parallelism on $P$ defined by $\ka$ we may associate with $\ka$ the
function 
$$\align
k:&P\to \h\otimes\La^2\h^*\\
k(u)(X,Y):&= K(\ze_X(u),\ze_Y(u))\text{ for }u\in P\text{ and } 
X,Y\in \h.
\endalign$$
We say that the Cartan connection $\ka$ has \idx{\it constant 
curvature} if this function $k$ is constant.

\head\totoc\nmb0{5}. Flat Cartan connections associated with a flat 
$G$-structure \endhead

\subhead{\nmb.{5.1}. $G$-structures}\endsubhead
By a $G$-structure on a smooth finite dimensional manifold $M$ we 
mean a principal fiber bundle $p:P\to M$ together with a 
representation $\rh: G\to GL(V)$ of the structure group
in a real vector space $V$ of dimension $\dim M$ and a 1-form $\si$ 
(called the \idx{\it soldering form}) on $M$ with values in the 
associated bundle $P[V,\rh]=P\x_G V$ which is fiber wise an 
isomorphism and identifies $T_xM$ with $P[V]_x$ for each $x\in M$. 
Then $\si$ corresponds uniquely to a $G$-equivariant 1-form 
$\th\in\Om^1_{\text{hor}}(P;V)^G$ which is strongly horizontal in the 
sense that its kernel is exactly the vertical bundle $VP$. The form 
$\th$ is called the \idx{\it displacement form} of the $G$-structure.

If the representation $\rh:G\to GL(V)$ is faithful so that 
$G\subset GL(V)$ is a linear Lie group, then a $G$-structure 
$(P,p,M,G,V,\th)$ is a subbundle of the linear frame bundle 
$GL(V;TM)$ of $M$. To see this recall the projection 
$q:P\x V\to P\x_GV=P[V]$ onto the associated bundle and the mapping 
$\ta=\ta^V:P\x_M P[V] \to V$ which is uniquely given by 
$q(u_x,\ta(u_x,v_x))= v_x$ and which satisfies $\ta(u_x,q(u_x,v))=v$ 
and $\ta(u_x.g,v_x)=\rh(g\i)\ta(u_x,v_x)$, see \cit!{14},~10.7. Then 
we have the smooth mapping over the identity on $M$,
$$\align
P &\to GL(V;P[V]) @>{\si\i}>> GL(V;TM),\tag1\\
u_x &\mapsto \ta(u_x,\quad)\i,
\endalign$$
which is $G$-equivariant and thus an embedding.

Let $p:P\to M$ and $p':P'\to M'$ be two $G$-structures with the same 
representation $\rh:G\to GL(V)$, with soldering 
forms $\si$, $\si'$, and with displacement forms $\th$, $\th'$, 
respectively. We are going to define the notion of an 
\idx{\it isomorphism of $G$-structures}. 

For $G\subset GL(V)$ an isomorphism of $G$-structures is a 
diffeomorphism $\ph:M\to M'$ such that the natural prolongation 
$GL(V;T\ph):GL(V;TM)\to GL(V;TM')$ to the frame bundles maps the 
subbundle $P\subset GL(V;TM)$ to the subbundle $P'\subset GL(V;TM')$.

In the general case an isomorphism of $G$-structures with the same 
representation $\rh:G\to GL(V)$ is an isomorphism of principal $G$-bundles
$$\CD
P @>\bar\ph>> P' \\
@VpVV     @VVp'V\\
M @>>\ph>  M
\endCD\tag2$$ 
whose induced isomorphism $\bar\ph\x_GId_{\rh(G)}$ of of principal 
$\rh(G)$-bundles coincides on $P\x_G\rh(G)\subset GL(V,TM)$ with the 
restriction of the natural prologation $GL(V;T\ph)$ of $\ph$ to the 
linear frame bundle, and
which preserves the displacement forms: 
$\bar\ph^*\th'=\th'\o T\bar\ph = \th:TP\to V$. 

A $G$-structure $(P,p,M,G,V,\th)$ is called \idx{\it flat} if it is 
locally (in a neighborhood of any point $x\in M$) isomorphic to the 
standard flat $G$-structure $\operatorname{pr_1}:V\x G\to V$ with 
displacement form $d\operatorname{pr_1}:T(V\x G)\to V$. 
Then the soldering form is just the identity
$\si=Id:TV=V\x V\to V\x V$, the linear frame bundle is 
$GL(V;TV)=V\x GL(V)$ and the associated $\rh(G)$-bundle is the 
subbundle $V\x\rh(G)\subset V\x GL(V)$. 

The standard 
examples of flat $G$-structures are foliations with structure group 
$$\pmatrix GL(p) & * \\ 0 & GL(n-p) \endpmatrix,$$
and symplectic structures.

Suppose that $G=\rh(G)\subset GL(V)$ and let us consider 
a local diffeomorphism defined near 0 and respecting 0 in $V$ (we 
write $\ph:V,0\to V,0$). It is a local automorphism of the 
standard flat $G$-structure $\operatorname{pr_1}:V\x G\to V$ with 
displacement form $d\operatorname{pr_1}:T(V\x G)\to V$ if and only if 
the following condition holds:
\roster
\item[3] $d\ph(x):V\to V$ is in $G\subset GL(V)$ for all $x$ in the 
       domain of $\ph$, 
\endroster
because only then its natural prolongation $GL(V;T\ph)$ to the linear 
frame bundle maps $G$-frames to $G$-frames.

\subhead\nmb.{5.2}. The infinite dimensional Lie group $GL^\infty(V)$
\endsubhead
Let $V$ be a real vector space of dimension $n$, and let 
$J^\infty(V,V)$ be the linear space of all infinite jets of smooth 
mappings $V\to V$, equipped with the initial topology with respect to 
all projections $J^\infty(V,V)\to J^k(V,V)$, which is a nuclear 
Fr\'echet space topology. 

We shall use the calculus of Fr\"olicher 
and Kriegl in infinite dimensions, see \cit!{9}, \cit!{15}, 
\cit!{16}, where a smooth mapping is one which maps smooth curves to 
smooth curves. On the spaces which we are going to use here the 
smooth mappings with values in finite dimensional spaces 
are just those which locally factor over some finite 
dimensional quotient like $J^k(V,V)$ and are smooth there.

Then we consider the closed linear subspace 
$\frak g\frak l_\infty(V)\subset J^\infty(V,V)$ of infinite jets of 
smooth mappings $V\to V$ which map the origin to the origin. Note 
that composition is defined on $\frak g\frak l_\infty(V)$ and is 
smooth, but is linear only in one (the left) component. Then we 
consider the open subset $GL_\infty(V)$ of all infinite jets of local 
diffeomorphisms of $V$, defined near and respecting 0. This is a 
smooth Lie group in the sense that composition and inversion are 
smooth. Its Lie algebra is $\frak g\frak l_\infty(V)$ which we may 
view as the full prolongation
$$
\frak g\frak l_\infty(V)=\frak g\frak l(V)\x \frak g\frak l^1(V) 
\x \frak g\frak l^2(V) \x \frak g\frak l^3(V)\x\dotso,\tag1
$$
where $\frak g\frak l^k(V)=S^kV^*\otimes V$ is the space of 
homogeneous polynomials $V\to V$ of order $k$. One may view
$\frak g\frak l_\infty(V)$ also as the vector space 
$\{j^\infty_0X: X\in\X(V), X(0)=0\}$ with the bracket 
$$[j^\infty_0X,j^\infty_0Y] = - j^\infty_0[X,Y]$$
and with the smooth (unique) exponential mapping 
$\exp:\frak g\frak l_\infty(V)\to GL_\infty(V)$ given by
$$\exp(j^\infty_0X) = j^\infty_0(\Fl^X_1),$$
where $\Fl^X_t$ is the flow of the vector field $X$ on $V$.
It is well known that $\exp:\frak g\frak l_\infty(V)\to GL_\infty(V)$ 
is not surjective onto any neighborhood of the identity, see 
\cit!{23}.

See \cit!{14}, section 13, for a detailed discussion of the finite jet 
groups $GL_k(V)$; the book \cit!{16} will contain a thorough 
discussion of $GL_\infty(V)$.

\subhead\nmb.{5.3}. The infinite prolongation of a linear Lie group 
$G$ and its Lie algebra\endsubhead
Let $G\subset GL(V)$ be a closed linear Lie group. We denote by 
$G_\infty\subset GL_\infty(V)$ the subgroup of all infinite jets 
$j^\infty_0\ph$ of local automorphisms $\ph$ of the standard flat 
$G$-structure $\operatorname{pr_1}:V\x G\to V$, defined near $0$ and 
respecting 0. Note that these $\ph$ are exactly the local diffeomorphisms 
$\ph:V,0\to V,0$ such that $d\ph(x)\in G\subset GL(V)$ for all 
$x\in V$ near 0, by the discussion in \nmb!{5.1}. Then $G_\infty$ is 
a group with multiplication and inversion
$$\gather
j^\infty_0\ph\o j^\infty_0\ps := j^\infty_0(\ph\o\ps),\\
(j^\infty_0\ph)\i = j^\infty_0(\ph\i),
\endgather$$
respectively. We will not address the question here in which sense 
$G_\infty$ is a Lie group. We continue just on a formal level.

The infinitesimal automorphisms respecting 0 of the standard flat 
$G$-structure on $V$ are then those local vector fields $X$ defined 
near 0 and vanishing at 0 in $V$ whose local flows $\Fl^X_t$ consist of 
automorphisms of the $G$-structure. 

\proclaim{Lemma}
The infinitesimal automorphisms are exactly the 
vector fields $X$ defined near 0 and vanishing near 0 such that $dX(x)\in \g\in L(V,V)$, where $\g$ is 
the Lie algebra of $G$.
\endproclaim

\demo{Proof}
Namely, $c(t) = d(\Fl^X_t)(x)$ is a curve in 
$G\subset GL(V)$ if and only if the following expression lies in 
$\g$:
$$\align
c'(t).c(t)\i &= \frac d{dt}(d(\Fl^X_t)(x)).d(\Fl^X_t)(x)\i
     = d(\frac d{dt}\Fl^X_t)(x).d(\Fl^X_t)(x)\i\\
&= d(X\o\Fl^X_t)(x).d(\Fl^X_t)(x)\i\\
&= d(X(\Fl^X_t)(x)).\qed
\endalign$$
\enddemo

We consider now the infinite jets 
$j^\infty_0X$ of all these infinitesimal automorphisms respecting 0.
These jets form a sub vector space 
$\g_\infty\subset \g\frak l_\infty(V)$ which we may view as the full 
prolongation 
$$
\frak g_\infty=\frak g\x \frak g^2 
\x \frak g^3 \x \frak g^4\x\dotso,\tag1
$$
where $\g^k=\g_\infty\cap S^kV^*\otimes V$ is the space of 
homogeneous polynomials $V\to V$ of order $k$ in $\g_\infty$. 
Then $\g_\infty$  is a Lie algebra with the bracket 
$$[j^\infty_0X,j^\infty_0Y] = - j^\infty_0[X,Y]$$
and with the smooth (unique) exponential mapping 
$\exp:\frak g_\infty\to G_\infty$ given by
$$\exp(j^\infty_0X) = j^\infty_0(\Fl^X_1),$$
where $\Fl^X_t$ is the flow of the vector field $X$ on $V$.
We expect that in general $\exp:\frak g_\infty\to G_\infty$ 
is not surjective onto any neighborhood of the identity.

Now we consider the Lie algebra of all infinitesimal automorphisms of 
the standard flat $G$-structure $\operatorname{pr}_1:V\x G\to V$, 
i\.e\. all local vector fields $X$ defined near 0 in $V$ such that 
the local flows $\Fl^X_t$ are automorphisms. As above one sees that 
these are the vector fields $X$ with $dX(x)\in\g\subset\g\frak l(V)$ 
for all $x$, without the restriction that they should vanish at 0.
Let $\frak a_\infty$ be the Lie algebra of all infinite jets 
$j^\infty_0X$ of such fields, again with bracket
$$
[j^\infty_0X,j^\infty_0Y] = - j^\infty_0[X,Y].
$$
By decomposing into monomials we have again $$
\frak a_\infty=V\x \frak g\x \frak g^2 \x \frak g^3 
\x \frak g^4\x\dotso = V\oplus \g_\infty.\tag2
$$
We have an adjoint representation 
$\Ad:G_\infty\to \operatorname{Aut}(\frak a_\infty)$ which is given by 
$$
\Ad(j^\infty_0\ph)j^\infty_0X = j^\infty_0(\ph^*X) = 
j^\infty_0(T\ph\i\o X\o \ph)
$$
In a formal sense we have also the left Maurer-Cartan form on 
$G_\infty$. First let us define the tangent bundle $TG_\infty$ as the 
set of all $(j^\infty_0\ph_0,j^\infty_0\frac d{dt}|_0\ph_t)$ where 
$\ph_t$ is a smooth curve of local automorphisms of the standard flat 
$G$-structure, respecting 0, smooth in the sense that 
$(t,x)\mapsto \ph_t(x)$ is smooth. Then we define the \idx{\it left} 
Maurer-Cartan form $\ka^l_{G_\infty}$ by
$$
\ka^l_{G_\infty}(j^\infty_0\ph_0,j^\infty_0\tfrac d{dt}|_0\ph_t) := 
j^\infty_0(T\ph_0\i\o\tfrac d{dt}|_0\ph_t).\tag3
$$

\proclaim{\nmb.{5.4}. Proposition}
Let $\th_0=d\operatorname{pr}_1:T(V\x G)\to V$ be the displacement 
form of the standard flat $G$-structure. Then 
$$\ka_0 := \th_0 \oplus \ka^l_{G_\infty}:T(V\x G_\infty) = 
TV\x TG_\infty \to V\oplus g_\infty = \frak a_\infty$$
is a flat Cartan connection on the manifold $V\x G_\infty$ with 
values in the Lie algebra $\frak a_\infty$.
\endproclaim

\demo{Proof}
Note first that the left Maurer-Cartan form $\ka^l_{G_\infty}$ given 
in \nmb!{5.3}.\thetag3 is really a trivialization of the tangent 
bundle $TG_\infty$ because of the lemma in \nmb!{5.3}, and we show 
that it satisfies the Maurer-Cartan equation:

Let $X$ (and later also $Y$) be a local vector field defined near 0 
and vanishing at 0, which is an infinitesimal automorphism of the 
standard flat $G$-structure, so that $j^\infty_0X$ is a typical element in 
$\g_\infty$. Let $j^\infty_0\ph\in G_\infty$ be a typical element, so 
$\ph:V,0\to V,0$ is a local automorphism. 
Then
$$
L_{j^\infty_0X}:G_\infty\ni j^\infty_0\ph\mapsto 
     T(j^\infty_0\ph)\o j^\infty_0X =j^\infty_0(T\ph\o X)
$$
is the left invariant vector field on $G_\infty$ generated by 
$j^\infty_0X$. We have 
$$
[L_{j^\infty_0X},L_{j^\infty_0Y}]
     =L_{[j^\infty_0X,j^\infty_0Y]}=L_{-j^\infty_0[X,Y]},
$$
and the Maurer-Cartan equation follows as usual:
$$\multline
(d\ka^l_{G_\infty})(L_{j^\infty_0X},L_{j^\infty_0Y}) 
  %   &= 
  %   L_{j^\infty_0X}(\ka^l_{G_\infty}(L_{j^\infty_0Y})) - 
  %   L_{j^\infty_0Y}(\ka^l_{G_\infty}(L_{j^\infty_0X})) -
  %   \ka^l_{G_\infty}([L_{j^\infty_0X},L_{j^\infty_0Y}])\\
= 0 - 0 - \ka^l_{G_\infty}(L_{[j^\infty_0X,j^\infty_0Y]}) 
     = -[j^\infty_0X,j^\infty_0Y] =\\
= -[\ka^l_{G_\infty}(L_{j^\infty_0X}),\ka^l_{G_\infty}(L_{j^\infty_0Y})]
     = -\frac12 [\ka^l_{G_\infty},\ka^l_{G_\infty}]^\wedge 
     (L_{j^\infty_0X},L_{j^\infty_0Y}).
\endmultline$$
One may now easily carry over to $\ka_0$ this result.
\qed\enddemo

\subhead\nmb.{5.5}. The infinite prolongation of a flat $G$-structure
\endsubhead
Let again $G\subset GL(V)$ be a linear Lie group, and let 
$(P,p,M,G,V,\th)$ be a flat $G$-structure. We denote by 
$p_\infty:P^\infty\to M$ the \idx{\it infinite prolongation} of this 
$G$-structure which is defined as follows:

The total space $P^\infty$ is the space of all infinite jets 
$j^\infty_0\ph$ of local isomorphisms $\ph:V,0\to M$ of the standard 
flat $G$-structure onto the given one. Then obviously the group 
$G_\infty$ acts freely from the right on $P^\infty$, and also 
transitive on the fiber. We have the mapping 
$$\gather
\ta^{P^\infty}:P^\infty\x_MP^\infty\to G_\infty\\
\ta^{P^\infty}(j^\infty_0\ph,j^\infty_0\ps) := j^\infty_0(\ph\i\o\ps)
\endgather$$
(see \cit!{14},~10.2) describing the principal $G_\infty$-bundle 
structure, which is locally isomorphic to the trivial bundle 
$V\x G_\infty$. The local isomorphisms $\ph:V,0\to M$ induce on 
$P^\infty$ a flat Cartan connection
$$\ka:TP^\infty\to\frak a_\infty = V\oplus \g_\infty$$
which locally is just given 
as the push forward via $j^\infty_0\ph$ of the canonical flat Cartan 
connection $\ka_0$ on $V\x G_\infty$. 

\head\totoc\nmb0{6}. The canonical Cartan connection for
a $G$-structure of first or second order \endhead

\subhead\nmb.{6.1} \endsubhead Let $G\subset GL(V)$ be 
a linear Lie group with Lie algebra $\g$.
We assume that the $G$-module $V\otimes \La^2V^*$ admits a decomposition 
$$V\otimes \La^2V^* = \de(\g\otimes V^*)\oplus \frak d\tag1$$ 
where $\de:V\otimes V^*\otimes V^*\to V\otimes\La^2V^*$ is the 
Spencer operator of alternation, and where $\frak d$ is a 
$G$-submodule.

We now recall the definition of the torsion function of a 
$G$-structure and the construction of its first prolongation in the 
sense of \cit!{22}. Let $p:P\to M$ be a $G$-structure on a manifold 
$M$ with a displacement form $\th:TP\to V$, see \nmb!{5.1}.
The 1-jet $j^1_xs$ of a local section $s:M\supset U\to P$ near 
$x\in M$ may be identified with its image $j^1_xs(T_xM)=H$, a 
horizontal linear subspace $H\subset T_{s(x)}P$. So the first jet 
bundle $J^1P\to P$ may be identified with the space of all horizontal 
linear subspaces in fibers of $TP\to P$.

Then the restriction $\th|H$ to such a horizontal space of the 
displacement form $\th$ is a linear isomorphism $\th|H:H\to V$ and we 
may use it to define the torsion function
$$\gather
t:J^1(P)\to V\otimes \La^2V^*\\
t(H)(v,w):= d\th((\th|H)\i(v),(\th|H)\i(w))
\endgather$$
We consider $P^1:= t\i(\frak d)$. It is a sub fiber bundle of $J^1(P)$ 
and the abelian vector group 
$G^1:=\g\otimes V^*\cap V\otimes S^2V^*\subset \operatorname{Hom}(V,\g)$
acts on $P^1$ freely by 
$g^1: P^1\ni H\mapsto g^1(H):= \{h + \ze^P_{g^1(\th(h))}(p(H)):h\in H\}$
where $\ze^P:\g\to \X(P)$ is the fundamental vector field mapping.
The orbits of $P^1$ under this $G^1$-action are fibers of the natural 
projection $p^1:P^1\to P$, hence $p^1:P^1\to P$ becomes a 
principal $G^1$-bundle.

Moreover, there exists a natural displacement form $\th^1$ on $P^1$. 
In order to define it we denote by $\Ph_H:T_uP\to V_uP$ the 
projection onto the vertical bundle $V_uP$ along the horizontal 
subspace $H\subset T_uP$. Then we have a well defined $\g$-valued 
$p^1$-horizontal 1-form $\om\in\Om^1(P^1;\g)^G$ given by 
$\ze^P_{\om_H(X)}(p^1(H)) = \Ph_H(T_H(p^1).X)$. It is part of the 
\idx{\it universal connection form} on the bundle of all connections 
$J^1(P)$.
The 1-form 
$$\th^1= \om + \th\o T(p^1): TP^1 \to \g \ltimes V$$
with values in the semidirect product $\g\ltimes V$ 
is the desired displacement form. It is equivariant with respect to 
the free action of the semidirect product $G\ltimes G^1$, where $G^1$ 
acts on $\g\x V$ by $(g^1,(X,v))\mapsto (X+g^1(v),v)$. 
So we have proved the main parts of 

\proclaim{\nmb.{6.2}. Lemma}
The fibration $p_1:P^1\to M$ is a principal fiber 
bundle with structure group $G\ltimes G^1$. 

The fibration $p^1:P^1\to P$ is a principal bundle with structure 
group $G^1$ and a $G^1$-structure on $P$ with the displacement form 
$\th^1$. Moreover the form $\th^1:TP^1\to \g\ltimes V$ is 
$G\ltimes G^1$-equivariant.
\endproclaim

\subhead\nmb.{6.3}. $G$-structures of type 1 \endsubhead
We assume now that the first prolongation $G^1$ of the group $G$ is 
trivial. Any $G$-structure with such a structure group $G$ is then 
called a \idx{\it $G$-structure of type 1}. In this case the 
projection $p^1:P^1\to P$ is a diffeomorphism and the displacement 
form $\th^1=\om + (p^1)^*\th$ may be identified with a 
$G$-equivariant 1-form on $P$ with values in $V_1:= \g \ltimes V$. 
This form is a Cartan connection in the principal bundle $p:P\to M$ 
of type $(\g\ltimes V)/\g$.

\proclaim{Proposition}
Let $G\subset GL(V)$ be a linear Lie group of type 1 satisfying 
condition \nmb!{6.1}.\thetag1. Then for any $G$-structure 
$(P,p,M,G,\th)$ the displacement form $\th^1$ of the first 
prolongation $p^1:P^1\to P$ defines a Cartan connection on $P$ of 
type $(\g\ltimes V)/\g$.
\endproclaim

\subhead\nmb.{6.4}. $G$-structures of type 2 \endsubhead
We say that $G\subset GL(V)$ is a \idx{\it linear Lie group of type 
2}, if the first prolongation $G^1\subset GL(V_1)$ is not trivial, but 
the second prolongation $(G^1)^1$ is trivial. So we have 
$$
\g\otimes V^* \cap V\otimes S^2V^* \ne 0, \quad
\g\otimes S^2V^* \cap V\otimes S^3V^* =0.
$$
We assume also that the condition \nmb!{6.1}.\thetag1  is satisfied.

Let $(P,p,M,G,\th)$ be a $G$-structure and let $p^1:P^1\to P$ be its 
first prolongation with the displacement form 
$\th^1:TP^1\to V_1=\g\ltimes \g$. We note that the 
$G^1$-submodule $\de(\g^1\otimes V_1^*)$ of $V_1\otimes \La^2V_1^*$ 
is not a direct summand. However, we may assume (at least when $G$ is 
reductive) that there exists a $G$-submodule $\frak d^1$ such that 
$$
V_1\otimes \La^2V_1^* = \de(\g^1\otimes V_1^*)\oplus \frak d^1.\tag1
$$
This is the case if there exist $G$-submodules $\frak d_1$, 
$\frak d_2$, $\frak d_3$ such that
$$\align
\g\otimes \La_2V_* &= R(\g)\oplus \frak d_1,\\
R(\g) &= \de(\g^1\otimes V^*)\oplus \frak d_2, \\
\g\otimes V^* &= \g^1\oplus \frak d_3,
\endalign$$
where $R(\g)$ is the space of curvature tensors of type $\g$, i\.e\. 
the space of closed (with respect to the Spencer  
differential) $\g$-valued 2-forms, and where $\frak d_2$ may be 
identified with the second Spencer cohomology space.

We denote by $t_1:J^1P^1\to V_1\otimes\La^2V_1^*$ the torsion 
function of the $G^1$-structure $p^1:P^1\to P$. The inverse image 
$P^2=t_1\i(\frak d^1)$ defines a submanifold of $J^1P^1$. The natural 
projection $p^2:P^2\to P^1$ is a diffeomorphism since by assumption 
the second prolongation $G^2$ of the group $G$ is trivial. In other 
words, we have a canonical field of horizontal subspaces in $TP^1$. 
Note that it is not a principal connection since it is not invariant 
under the group $G^1$. Using this field of horizontal subspaces in 
$TP^1$ we may extend the canonical vertical parallelism 
$VP^1\to \g\ltimes\g^1$ of the principal $G\ltimes G^1$-bundle 
$p_1:P^1\to M$ to a $G$-quivariant 1-form $\om_1$ on $P^1$. The 
$\g$-component of $\om_1$ is the $\g$-valued 1-form $\om$ from 
\nmb!{6.1}. 

The form 
$$
\th^2 = (p^1)^*\th + \om_1 : TP^1\to V \x \g \x \g^1 = \frak a_\infty
$$
is non degenerate, $G$-equivariant, and it prolongs the vertical 
parallelism $VP^1\to \g\ltimes \g^1$ of the principal 
$G\ltimes G^1$-bundle $p_1:P^1\to M$. Hence it is a Cartan connection 
of type $\frak a_\infty/(\g\ltimes \g^1)$, where 
$\frak a_\infty= V\x \g \x \g^1$ is the full prolongation of the Lie 
algebra $\g\subset \g\frak l(V)$.

Summarizing we have

\proclaim{\nmb.{6.5}. Proposition}
Let $p:P\to M$ be a $G$-structure of type 2 satisfying the conditions 
\nmb!{6.1}.\thetag1 and \nmb!{6.4}.\thetag1. 

Then on the total space $P^1$ of the first prolongation 
$p^1:P^1\to P$ of the bundle $p:P\to M$ there exists a canonically 
defined Cartan connection of type $\frak a_\infty/\g$, where 
$\frak a_\infty = V \x \g \x \g^1$ is the full prolongation of the 
Lie algebra $\g$ of $G$.
\endproclaim

\enddocument